\documentclass[11pt,onecolumn,draftcls]{IEEEtran}

\usepackage{ifpdf}

\usepackage{cite}

\ifCLASSINFOpdf
 \usepackage[pdftex]{graphicx}
 \DeclareGraphicsExtensions{.pdf,.jpeg,.png}
\else
 \usepackage[dvips]{graphicx}
 \DeclareGraphicsExtensions{.eps}
\fi

\usepackage[cmex10]{amsmath}

\usepackage{algorithmic}
\usepackage{array}
\usepackage{mdwmath}
\usepackage{mdwtab}
\usepackage{eqparbox}

\usepackage[tight,footnotesize]{subfigure}

\usepackage{stfloats}
\usepackage{url}

\usepackage{amssymb}
\usepackage{bm}
\usepackage{overpic}
\usepackage{color}
\usepackage{dsfont}
\usepackage{amsthm}

\newtheorem{theorem}{Theorem}
\newtheorem{lemma}{Lemma}

\newtheorem{assumption}{Assumption}
\newtheorem{example}{Example}

\def\E{\mathbb{E}}
\def\mc{\mathcal}

\title{Diffusion Adaptation Strategies for Distributed Optimization and Learning over Networks}
\begin{document}

\author{Jianshu~Chen,~\IEEEmembership{Student~Member,~IEEE,}%
        ~and~Ali~H.~Sayed,~\IEEEmembership{Fellow,~IEEE}
\thanks{Manuscript received October 30, 2011; revised March 15, 2012. 
This work was supported in part by NSF grants CCF-1011918 and CCF-0942936.
Preliminary results related to this work are reported in the conference presentations%
\cite{chen2011diffusionOpt} and \cite{chen2011MSE}.}%
\thanks{The authors are with Department of Electrical Engineering,
University of California, Los Angeles, CA 90095. Email: \{jshchen, sayed\}@ee.ucla.edu.
}
}

\maketitle

\begin{abstract}
We propose an adaptive diffusion mechanism to optimize a global cost function
in a distributed manner over a network of nodes. The cost function is assumed to
consist of a collection of individual components. Diffusion adaptation allows the
nodes to cooperate and diffuse information in real-time; it also helps alleviate the effects of
stochastic gradient noise and measurement noise through a continuous learning process.
We analyze the mean-square-error performance
of the algorithm in some detail, including its transient and steady-state behavior. We also apply the diffusion algorithm
to two problems: distributed estimation
with sparse parameters and
distributed localization.
Compared to well-studied incremental methods, diffusion
methods do not require the use of a cyclic path over the nodes and are
robust to node and link failure.
 Diffusion methods also endow networks with adaptation abilities that enable the individual nodes to continue learning even when the cost function changes with time. Examples involving such dynamic cost functions with moving targets are common in the context of biological networks.
\end{abstract}

 \ifCLASSOPTIONpeerreview
 \begin{center} \bfseries EDICS Category: SEN--DIST, SEN--COLB, ASD--ANAL \end{center}
 \fi


\section{Introduction}
\label{sec:intro}
\IEEEPARstart{W}{e} consider the problem of optimizing  a global  cost function in a distributed manner.
The cost function is assumed to consist of the sum of individual components, and
spatially distributed nodes are used to seek the common minimizer (or maximizer)
through local interactions.
{\color{black}
Such problems abound in the context of biological networks, where
agents collaborate with each other via local interactions
for a common objective, such as locating food sources or evading predators\cite{tu2011fish}.
Similar problems are common in distributed resource allocation applications and in
online machine learning procedures. In the latter case, data that are generated by the same underlying distribution are processed in a distributed manner over a network of learners in order to recover the model parameters
(e.g., \cite{dekel2011optimal,zaid2011collaborativeGMM}).
}%

There are already a few of useful techniques for the solution of
optimization problems in a distributed manner
\cite{bertsekas1997new, nedic2001incremental, rabbat2005quantized,
lopes2007incremental,bertsekasparallel,tsitsiklis1984convergence,tsitsiklis1986distributed,barbarossa2007bio,
nedic2009bookchapter,
nedic2009distributed,
schizas2008consensus1,kar2008sensor,kar2011converegence,
dimakis2010gossip,olfati2004consensus,
aysal2009broadcast,sardellitti2010fast,xiao2005scheme,eksin2011asilomar}.
Most notable among these
methods is the incremental approach
\cite{bertsekas1997new, nedic2001incremental, rabbat2005quantized,
lopes2007incremental}
and the consensus approach
\cite{bertsekasparallel,tsitsiklis1984convergence,tsitsiklis1986distributed,barbarossa2007bio,
nedic2009bookchapter,
nedic2009distributed,
schizas2008consensus1,kar2008sensor,kar2011converegence,
dimakis2010gossip,olfati2004consensus,
aysal2009broadcast,sardellitti2010fast,xiao2005scheme}.
In the incremental approach, a cyclic path is defined over the nodes and data are processed in a cyclic manner through
the network until optimization is achieved. However, determining a cyclic path that covers all nodes
is known to be an NP-hard problem\cite{karp1972reducibility}
and, in addition, cyclic trajectories are prone to link and node failures.
When any of the edges along the path fails, the sharing of data through the cyclic trajectory is interrupted and the algorithm stops performing.
In the consensus approach, vanishing step-sizes are used to ensure that nodes reach consensus and converge to the same optimizer in steady-state. However, in time-varying environments, diminishing step-sizes prevent the network from continuous learning and optimization; when the step-sizes die out, the network stops learning. In earlier publications
\cite{lopesdistributed,lopes2007diffusion,Sayed07,lopes2008diffusion,
cattivelli2008diffusion,Cattivelli10,cattivelli2007diffusionRLS,cattivelli2008TSPdiffusionRLS,cattivelli2010TACdiffusionKalman,
takahashi2010diffusion}, and motivated by our work on adaptation and learning over networks,
we introduced the concept of diffusion adaptation and showed how this technique can be used to
solve global minimum mean-square-error estimation problems efficiently \emph{both}
in real-time \emph{and} in a distributed manner.
In the diffusion approach, information
is processed locally and \emph{simultaneously} at all nodes and the processed data are
diffused through a real-time sharing mechanism that ripples through the network continuously.
Diffusion adaptation was applied to model complex patterns of behavior encountered
in biological networks, such as bird flight formations \cite{cattivelli2011modeling} and
fish schooling \cite{tu2011fish}.
Diffusion adaptation was also applied to solve dynamic resource allocation problems in cognitive radios\cite{di2011bio}, to perform robust system identification\cite{chouvardas2011adaptive},
and to implement distributed online learning
in pattern recognition applications\cite{zaid2011collaborativeGMM}.

This paper generalizes the diffusive learning process and applies it to the distributed
optimization of a wide class of cost functions.
The diffusion approach will be shown to alleviate the effect of gradient noise on convergence.
Most other studies on distributed optimization tend to focus on the almost-sure convergence
of the algorithms under diminishing step-size conditions\cite{bertsekas1997new,nedic2001incremental,
ram2010distributed,bianchi2011convergence,bertsekas2010incremental,borkar2000ode,srivastava2011distributed}, or
on convergence under deterministic conditions on the data
\cite{bertsekas1997new,nedic2001incremental,rabbat2005quantized,nedic2009distributed}.
In this article we instead examine the distributed algorithms from a mean-square-error perspective at \emph{constant} step-sizes.
This is because constant step-sizes are necessary for continuous adaptation, learning, and tracking, which in turn  enable the resulting algorithms to perform well even under data that exhibit statistical variations, measurement
noise, and gradient noise.

{\color{black}
This paper is organized as follows. In Sec. \ref{Sec:ProblemFormulation}, we  introduce the global cost function and approximate it by a distributed optimization problem through the use of a second-order Taylor series expansion. In Sec. \ref{Sec:Diffusion}, we show that optimizing the localized alternative
cost at each node $k$ leads naturally to diffusion adaptation strategies.
In Sec. \ref{Sec:ConvergenceAnalysis}, we analyze the mean-square performance of the diffusion
algorithms under statistical perturbations when stochastic gradients are used.
In Sec. \ref{Sec:Simulation}, we apply the diffusion algorithms to two application
problems: sparse distributed estimation and distributed
localization. Finally, in Sec. \ref{Sec:Conclusion}, we conclude the paper.
}%

\noindent
{\bf Notation}. Throughout the paper, all vectors are column vectors except for the regressors $\{\bm{u}_{k,i}\}$, which are taken to be row vectors for simplicity of notation.
We use boldface letters to denote random quantities (such as $\bm{u}_{k,i}$)
and regular font letters to denote their realizations or deterministic variables (such as $u_{k,i}$).
We write $\E$ to denote the expectation operator. We use $\mathrm{diag}\{x_1,\ldots,x_N\}$ to denote a diagonal matrix consisting of diagonal entries
$x_1,\ldots,x_N$,
and use $\mathrm{col}\{x_1,\ldots,x_N\}$ to denote a column vector formed by stacking $x_1,\ldots,x_N$
on top of each other. For symmetric matrices $X$ and $Y$, the notation $X \le Y$ denotes
$Y-X \ge 0$, namely, that the matrix difference $Y-X$ is positive semi-definite.

\section{Problem Formulation}
\label{Sec:ProblemFormulation}
The objective is to determine, in a collaborative and distributed manner,
the $M\!\times\! 1$ column vector  $w^o$ that minimizes a global cost
of the form:
	\begin{align}
		\label{Equ:ProblemFormulation:J_glob_original}
		\boxed{
			J^{\mathrm{glob}}(w)	=	\sum_{l=1}^N J_l(w)			
		}
	\end{align}
where $J_l(w)$, $l=1,2,\ldots,N$, are individual real-valued functions,
defined over $w \in \mathds{R}^M$ and assumed to be
differentiable and strictly convex.
Then, $J^{\mathrm{glob}}(w)$ in \eqref{Equ:ProblemFormulation:J_glob_original}
is also strictly convex so that the minimizer $w^o$ is unique\cite{poliak1987introduction}.
In this article we study the important case where the component functions $\{J_l(w)\}$
are minimized at the \emph{same} $w^o$. This case is common in practice; situations abound
where nodes in a network need to work cooperatively to attain a common objective
(such as tracking a target, locating the source of chemical leak, estimating a physical model,
or identifying a statistical distribution). This scenario is also frequent  in the context of biological networks.
For example, during the foraging behavior of an animal group, each agent in the group
is interested in determining the
\emph{same} vector $w^o$ that corresponds to the  location of the food source
or the location of the predator \cite{tu2011fish}.
This scenario is equally common in online distributed machine learning problems, where data samples
are often generated from the same underlying distribution and they are processed in a distributed manner by different nodes
(e.g., \cite{dekel2011optimal,zaid2011collaborativeGMM}).
The case where the $\{J_l(w)\}$ have different individual minimizers is studied in
\cite{chen2012ssp}; this situation is more challenging to study. Nevertheless,
it is shown in \cite{chen2012ssp} that the same diffusion strategies
\eqref{Equ:DiffusionAdaptation:ATC0}--\eqref{Equ:DiffusionAdaptation:CTA0}
of this paper are still applicable and
nodes would converge instead to a Pareto-optimal solution.

Our strategy to optimize the global cost $J^{\mathrm{glob}}(w)$ in a distributed manner
is based on three steps. First, using a second-order Taylor series expansion, we argue
that $J^{\mathrm{glob}}(w)$  can be approximated by an alternative localized cost that is amenable to
distributed optimization --- see \eqref{Equ:DiffusionAdaptation:J_glob_prime_prime_final}.
Second, each individual node optimizes this alternative cost via a
steepest-descent procedure that relies solely on interactions within the neighborhood of the node. Finally, the local estimates for $w^o$ are spatially combined by each node and the procedure repeats itself in real-time.

To motivate the approach, we start by introducing a set of nonnegative coefficients $\{c_{l,k}\}$ that satisfy:
	\begin{align}
		\label{Equ:ProblemFormulation:C_Condition}
		\boxed{
			 \displaystyle
			\sum_{k=1}^N c_{l,k} = 1,\quad
			c_{l,k}=0~\mathrm{if}~l \notin \mathcal{N}_k,		\quad  l=1,2,\ldots,N	
		}
	\end{align}
where $\mathcal{N}_k$ denotes the neighborhood of node $k$ (including node $k$ itself); the neighbors
of node $k$ consist of all nodes with which node $k$ can share information.
Each $c_{l,k}$ represents a weight value that node $k$ assigns to information arriving from its neighbor $l$.
Condition \eqref{Equ:ProblemFormulation:C_Condition} states that the sum of all weights leaving each
node $l$ should be one.
Using the coefficients $\{c_{l,k}\}$, we can express $J^{\mathrm{glob}}(w)$ from \eqref{Equ:ProblemFormulation:J_glob_original} as
	\begin{align}
		\label{Equ:ProblemFormulation:J_glob}
		J^{\mathrm{glob}}(w)  &=
                                        \displaystyle J_k^{\mathrm{loc}}(w) + \sum_{l \neq k}^N
                                        J_l^{\mathrm{loc}}(w)	
	\end{align}
where
	\begin{align}
		\label{Equ:ProblemFormulation:J_k_loc0}
		J_k^{\mathrm{loc}}(w)	\triangleq	\sum_{l \in \mathcal{N}_k} c_{l,k} J_l(w)	
	\end{align}
In other words, for each node $k$, we are introducing a new local cost function, $J_k^{\mathrm{loc}}(w)$,
which corresponds to a weighted combination of the costs of its neighbors.
Since the $\{c_{l,k}\}$ are all nonnegative and each $J_l(w)$ is convex, then $J_k^{\mathrm{loc}}(w)$ is also a
convex function (actually, the $J_k^{\mathrm{loc}}(w)$ will be guaranteed to be strongly convex in our
treatment in view of Assumption \ref{Assumption:Hessian}  further ahead).

Now, each $J_l^{\mathrm{loc}}(w)$ in the second term of \eqref{Equ:ProblemFormulation:J_glob}
can be approximated
via a second-order Taylor series expansion as:
	\begin{align}
		J_l^{\mathrm{loc}}(w)	
		\label{Equ:ProblemFormulation:J_k_loc}
							\approx\;
								&	J_l^{\mathrm{loc}}(w^o)
								+
								\|w-w^o\|_{\Gamma_l}^2
	\end{align}
where $\Gamma_l\!=\!\frac{1}{2}\nabla_w^2 J_l^{\mathrm{loc}}(w^o)$ is the
(scaled) Hessian matrix  relative to $w$ and evaluated at $w\!=\!w^o$,
and the notation
$\|a\|_\Sigma^2$ denotes $a^T \Sigma a$ for any weighting matrix $\Sigma$.
The analysis in the subsequent sections will show that the
second-order approximation \eqref{Equ:ProblemFormulation:J_k_loc}
is sufficient to ensure mean-square convergence
of the resulting diffusion algorithm.  Now, substituting
\eqref{Equ:ProblemFormulation:J_k_loc} into the right-hand side of
\eqref{Equ:ProblemFormulation:J_glob}
gives:
	\begin{align}
		\label{Equ:ProblemFormulation:J_glob_approx}
		J^{\mathrm{glob}}(w)	\approx	
							 \displaystyle
							J_k^{\mathrm{loc}}(w)
							\!+\!
							\sum_{l \neq k} \|w\!-\!w^o\|_{\Gamma_l}^2
							\!+\!
							\sum_{l \neq k} J_l^{\mathrm{loc}} (w^o)
	\end{align}
The last term in the above expression does not depend on the unknown $w$.
Therefore, we can ignore it so that optimizing $J^{\mathrm{glob}}(w)$
is approximately equivalent to optimizing the following alternative cost:
	\begin{align}
		\label{Equ:ProblemFormulation:J_glob_prime}
		J^{\mathrm{glob}'}(w)&\triangleq
							 \displaystyle
							J_k^{\mathrm{loc}}(w)
							+
							\sum_{l \neq k} \|w-w^o\|_{\Gamma_l}^2
	\end{align}

\section{Iterative Diffusion Solution}
\label{Sec:Diffusion}
Expression \eqref{Equ:ProblemFormulation:J_glob_prime} relates the original global cost
\eqref{Equ:ProblemFormulation:J_glob_original}
to the newly-defined local cost function $J_{k}^{\mathrm{loc}}(w)$. The relation is through the
second term on the right-hand side of \eqref{Equ:ProblemFormulation:J_glob_prime},
which corresponds to a sum of quadratic terms involving
the minimizer $w^o$. Obviously, $w^o$
is not available at node $k$ since the nodes wish to estimate $w^o$.
Likewise, not all Hessian matrices $\Gamma_l$ are
available to node $k$. Nevertheless, expression \eqref{Equ:ProblemFormulation:J_glob_prime}
suggests a useful approximation that leads to a powerful distributed solution, as we proceed to explain.

\begin{figure}[t]
    \centering
    \includegraphics[width=0.36\textwidth]{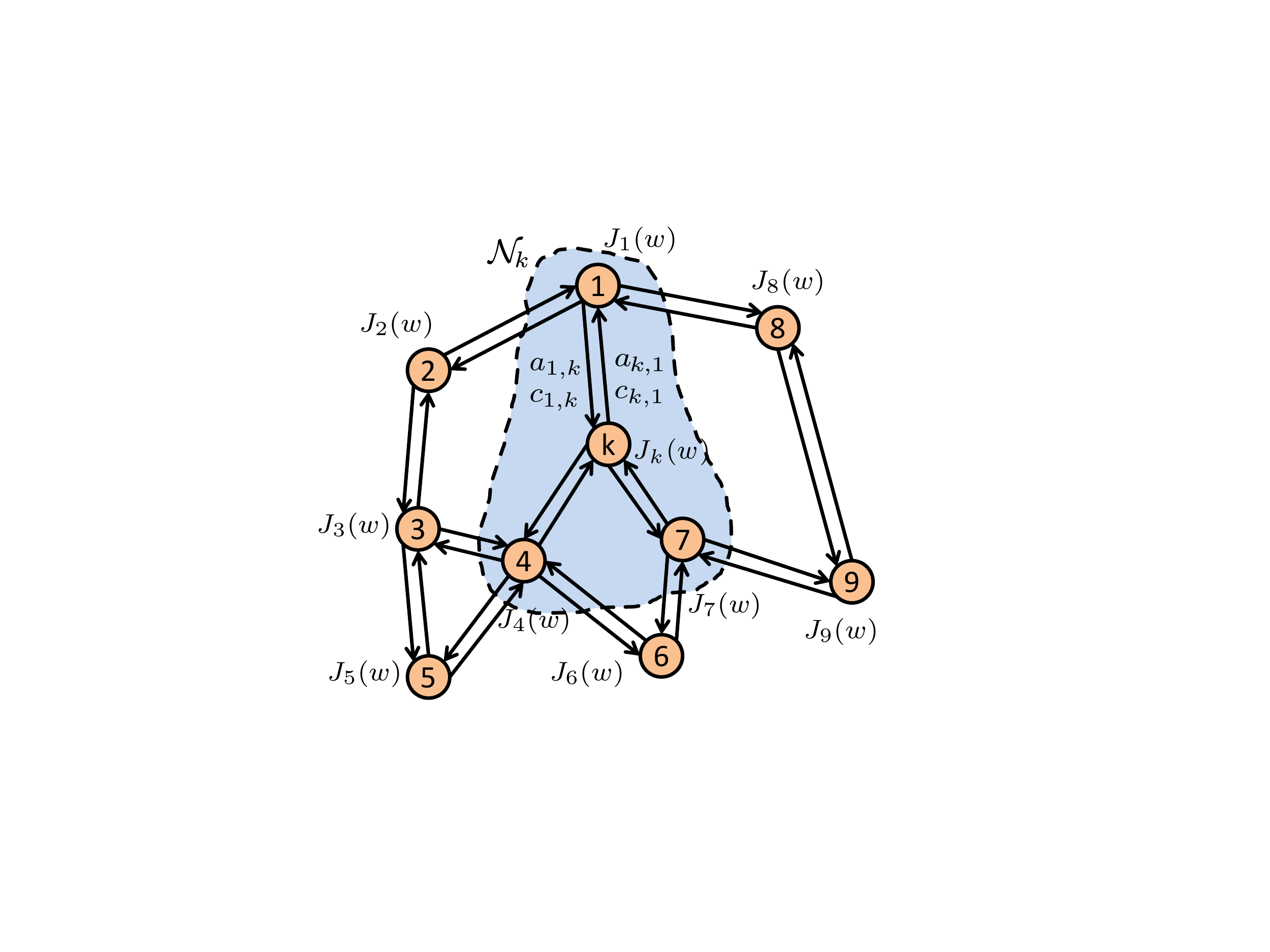}
    \caption{A network with $N$ nodes; a cost function $J_k(w)$ is associated with each node $k$. The set of neighbors of node $k$ is denoted by ${\cal N}_k$; this set consists of all nodes with which node $k$ can share information.}
    \label{Fig:Network}
\end{figure}

Our first step is to replace the global cost $J^{\mathrm{glob}'}(w)$ by a
reasonable \emph{localized} approximation for
it at every node $k$. Thus, initially we limit the summation on the right-hand side of \eqref{Equ:ProblemFormulation:J_glob_prime} to the neighbors of node $k$
 and introduce the cost function:
	\begin{align}
	 	\label{Equ:DiffusionAdaptation:J_glob_prime_local}
		J_k^{\mathrm{glob}'}(w) \triangleq  \displaystyle J_k^{\mathrm{loc}}(w)
							+
							\sum_{l\in{\mathcal{N}}_k\backslash\{k\}} \|w-w^o\|^2_{\Gamma_l}
	\end{align}
Compared with \eqref{Equ:ProblemFormulation:J_glob_prime},
the last term in \eqref{Equ:DiffusionAdaptation:J_glob_prime_local}
involves only quantities that are available in the neighborhood of node $k$. The argument involving steps
\eqref{Equ:ProblemFormulation:J_k_loc}--\eqref{Equ:DiffusionAdaptation:J_glob_prime_local}
therefore shows us one way by which we can adjust the earlier local cost function $J_k^\mathrm{loc}(w)$
defined in \eqref{Equ:ProblemFormulation:J_k_loc0} by adding to it the last term that appears in
\eqref{Equ:DiffusionAdaptation:J_glob_prime_local}.
{\color{black}
Doing so, we end up replacing
$J_k^\mathrm{loc}(w)$ by  $J_k^{\mathrm{glob}'}(w)$,
and this new localized cost function preserves the second term in \eqref{Equ:ProblemFormulation:J_glob}
up to a second-order approximation. This correction will help lead
to a diffusion step (see \eqref{Equ:DiffusionAdaptation:Combination_intermediate}--%
\eqref{Equ:DiffusionAdaptation:Combination_intermediate2}).
}

Now, observe that the cost in \eqref{Equ:DiffusionAdaptation:J_glob_prime_local}
includes the quantities $\{\Gamma_l\}$, which belong to the neighbors of node $k$.
These quantities may or may not be available. If they are known, then we can proceed with \eqref{Equ:DiffusionAdaptation:J_glob_prime_local} and rely
on the use of the Hessian matrices $\Gamma_l$
in the subsequent development. Nevertheless, the more interesting situation in practice is when these Hessian matrices are not known beforehand (especially since they depend on the unknown $w^o$). For this reason, in this article, we
approximate each $\Gamma_l$ in \eqref{Equ:DiffusionAdaptation:J_glob_prime_local}
by a multiple of the identity matrix,
say,
    \begin{align}
        \label{Equ:DiffusionAdaptation:Gamma_l_approx}
        \Gamma_l \approx b_{l,k} I_M
    \end{align}
for some nonnegative coefficients $\{b_{l,k}\}$; observe that we are allowing the
coefficient $b_{l,k}$ to vary with the node index $k$. Such approximations
are common in stochastic approximation theory and help reduce the complexity of the resulting algorithms  --- see \cite[pp.20--28]{poliak1987introduction} and \cite[pp.142--147]{Sayed08}.
Approximation \eqref{Equ:DiffusionAdaptation:Gamma_l_approx} is reasonable since, in view of the Rayleigh-Ritz characterization of eigenvalues \cite{golub1996matrix}, we can always bound the weighted
squared norm $\|w-w^o\|^2_{\Gamma_{l}}$ by the unweighted squared norm as follows
    \begin{align}
        \lambda_{\min}(\Gamma_{l})\cdot\|w\!-\!w^o\|^2\leq
        \|w\!-\!w^o\|^2_{\Gamma_{l}}\leq
        \lambda_{\max}(\Gamma_{l})\cdot \|w\!-\!w^o\|^2
        \nonumber
    \end{align}
Thus, we replace \eqref{Equ:DiffusionAdaptation:J_glob_prime_local}
by
	\begin{align}
		\label{Equ:DiffusionAdaptation:J_glob_prime_prime}
		J_k^{\mathrm{glob}''}(w) 	\triangleq	
								\displaystyle J_k^{\mathrm{loc}}(w)
								\;+\;
								\sum_{l\in{\mathcal N}_k\backslash\{k\}} b_{l,k} \|w-w^o\|^2
	\end{align}
As the derivation will show, we do not need to worry at this stage about how  the scalars $\{b_{l,k}\}$ are selected;
they will be embedded into other combination weights that the designer selects.
If we replace $J_k^\mathrm{loc}(w)$ by its definition \eqref{Equ:ProblemFormulation:J_k_loc0}, we can rewrite
\eqref{Equ:DiffusionAdaptation:J_glob_prime_prime} as
	\begin{align}
	\label{Equ:DiffusionAdaptation:J_glob_prime_prime_final}
		\boxed{
		J_k^{\mathrm{glob}''}\!(w)	= 	 \displaystyle
								\sum_{l\in{\cal N}_k} \! c_{l,k}J_l(w)+\!\!
								\sum_{l\in{\cal N}_k\backslash\{k\}}
								b_{l,k}\|w\!\!-\!\!w^o\|^2				
		}
	\end{align}
Observe that cost \eqref{Equ:DiffusionAdaptation:J_glob_prime_prime_final}
is different for different nodes; this is because the choices of the weighting scalars $\{c_{l,k},b_{l,k}\}$
vary across nodes $k$; moreover, the neighborhoods vary with $k$. Nevertheless, these localized
cost functions now constitute the important starting point for the development of diffusion
strategies for the online and distributed optimization of \eqref{Equ:ProblemFormulation:J_glob_original}.

Each node $k$ can apply a steepest-descent iteration to
minimize $J_k^{\mathrm{glob}''}(w)$ by moving along the negative direction of the gradient (column)
vector of the cost function, namely,
	\begin{align}
		w_{k,i}	=\;&	 \displaystyle w_{k,i-1}
                            - \mu_k \sum_{l \in \mathcal{N}_k} c_{l,k} \nabla_w J_l(w_{k,i-1})
                            \nonumber\\
		\label{Equ:DiffusionAdaptation:GradientDescent}
                \;&
					-
					 \displaystyle\mu_k \sum_{l \in \mathcal{N}_k\backslash \{k\}}
					 2b_{l,k}(w_{k,i-1}-w^o),
								\qquad i \ge 0				
	\end{align}	
where $w_{k,i}$ denotes the estimate for $w^o$ at node $k$ at time $i$, and $\mu_k$ denotes a small
\emph{constant}  positive
step-size parameter.
{\color{black}
While vanishing step-sizes, such as $\mu_{k}(i)=1/i$, can be used in
\eqref{Equ:DiffusionAdaptation:GradientDescent}, we consider in this paper the case of
constant step-sizes. This is because we are interested in distributed strategies
that are able to continue adapting and learning. An important question to address
therefore is how close each of the $w_{k,i}$ gets to the optimal solution
$w^o$; we answer this question later in the paper by means of a mean-square-error
convergence analysis (see expression
\eqref{Equ:ConvergenceAnalysis:MSD_Vanish_SmallStepSize_final}).
It will be seen then that the mean-square-error (MSE) of the algorithm will
be of the order of the step-size; hence, sufficiently small step-sizes will lead to
sufficiently small MSEs.
}

Expression \eqref{Equ:DiffusionAdaptation:GradientDescent}
adds two correction terms to the previous estimate, $w_{k,i-1}$, in order to update it to $w_{k,i}$.
The correction terms can be added one at a time in a succession of two steps, for example, as:
	\begin{align}
		\label{Equ:DiffusionAdaptation:Adaptation_intermedieate}
		\psi_{k,i}	&=	 \displaystyle w_{k,i-1} - \mu_k \sum_{l \in \mathcal{N}_k} c_{l,k} \nabla_w J_l(w_{k,i-1})	\\
		\label{Equ:DiffusionAdaptation:Combination_intermediate}
		w_{k,i}	&=	 \displaystyle \psi_{k,i}
					- \mu_k \sum_{l \in \mathcal{N}_k\backslash\{k\}} 2b_{l,k} (w_{k,i-1}-w^o)
	\end{align}
Step \eqref{Equ:DiffusionAdaptation:Adaptation_intermedieate}
updates $w_{k,i-1}$ to an intermediate value $\psi_{k,i}$ by using a \emph{combination} of local gradient vectors.
Step \eqref{Equ:DiffusionAdaptation:Combination_intermediate}
further updates $\psi_{k,i}$ to $w_{k,i}$  by using a \emph{combination} of local estimates.
However, two issues arise while examining
\eqref{Equ:DiffusionAdaptation:Combination_intermediate}:
	\begin{enumerate}
		\item[(a)]
		First, iteration \eqref{Equ:DiffusionAdaptation:Combination_intermediate}
		requires knowledge of the optimizer $w^o$.
        However, all nodes are running similar updates to estimate the $w^o$. By the time node $k$ wishes to apply \eqref{Equ:DiffusionAdaptation:Combination_intermediate}, each of its neighbors would have performed its own update similar to \eqref{Equ:DiffusionAdaptation:Adaptation_intermedieate}
        and would have available their intermediate estimates, $\{\psi_{l,i}\}$. Therefore,
		we replace $w^o$ in \eqref{Equ:DiffusionAdaptation:Combination_intermediate} by
		$\psi_{l,i}$.  This step helps diffuse information over the network
        and brings into node $k$ information that exists beyond its immediate neighborhood;
		this is because each $\psi_{l,i}$ is influenced by data from the neighbors of node $l$.
        {\color{black}
		We observe that this diffusive term arises from the quadratic approximation
		\eqref{Equ:ProblemFormulation:J_k_loc} we have made
		to the second term in \eqref{Equ:ProblemFormulation:J_glob}.
        }
		
		\item[(b)]
		Second, the intermediate value $\psi_{k,i}$ in
        \eqref{Equ:DiffusionAdaptation:Adaptation_intermedieate}
        is generally a better estimate for $w^o$
		than $w_{k,i-1}$ since it is obtained by incorporating information from the neighbors
		through \eqref{Equ:DiffusionAdaptation:Adaptation_intermedieate}.
		Therefore, we further replace $w_{k,i-1}$ in \eqref{Equ:DiffusionAdaptation:Combination_intermediate}
		by $\psi_{k,i}$. This step is reminiscent of
		incremental-type approaches to optimization, which have been widely
		studied in the literature
		\cite{bertsekas1997new,nedic2001incremental,
		rabbat2005quantized, lopes2007incremental}.
	\end{enumerate}
Performing the substitutions described in items (a) and (b) into
\eqref{Equ:DiffusionAdaptation:Combination_intermediate}, we obtain:
	\begin{align}
		\label{Equ:DiffusionAdaptation:Combination_intermediate2}
		w_{k,i}	&=	 \displaystyle \psi_{k,i}
					- \mu_k \sum_{l \in \mathcal{N}_k\backslash\{k\}} 2b_{l,k} (\psi_{k,i}-\psi_{l,i})
	\end{align}
Now introduce the coefficients
	\begin{align}
        \label{Equ:DiffusionAdaptation:a_lk_def}
		a_{l,k} \triangleq 2 \mu_k b_{l,k}	\quad (l \! \neq \! k),
		\quad
		a_{k,k} \triangleq 1 \!-\!  \mu_k \!\! \sum_{l \in \mathcal{N}_k\backslash\{k\}} 2b_{l,k}
	\end{align}
Note that the $\{a_{l,k}\}$ are nonnegative for $l\neq k$ and $a_{k,k}\geq 0$ for sufficiently small step-sizes. Moreover, the coefficients $\{a_{l,k}\}$ satisfy
    \begin{align}
		\label{Equ:AdaptiveDiffusion:Condition_a}
		 \displaystyle
		\sum_{l=1}^N a_{l,k} = 1, \quad a_{l,k} = 0~\mathrm{if}~l \notin \mathcal{N}_{k}		
	\end{align}
Using \eqref{Equ:DiffusionAdaptation:a_lk_def} in
\eqref{Equ:DiffusionAdaptation:Combination_intermediate2},
we arrive at the following Adapt-then-Combine (ATC)
diffusion strategy (whose structure is the same as the ATC algorithm originally proposed in
\cite{lopes2008diffusion,
cattivelli2008diffusion,Cattivelli10}
for mean-square-error estimation):
	\begin{align}
		\boxed{
				\label{Equ:DiffusionAdaptation:ATC0}
				\begin{array}{l}
					\psi_{k,i}	=	 \displaystyle
								w_{k,i-1} - \mu_k \sum_{l \in \mathcal{N}_k} c_{l,k}
								\nabla_w J_l(w_{k,i-1})	\\
					w_{k,i}	=	\displaystyle \sum_{l \in \mathcal{N}_k} a_{l,k} \psi_{l,i}
				\end{array}
			}
	\end{align}
To run algorithm \eqref{Equ:DiffusionAdaptation:ATC0},
we only need to select combination coefficients $\{a_{l,k},c_{l,k}\}$ satisfying
\eqref{Equ:ProblemFormulation:C_Condition}
and \eqref{Equ:AdaptiveDiffusion:Condition_a}, respectively; there is no need to worry about the
intermediate coefficients $\{b_{l,k}\}$ any more, since they have been blended into the $\{a_{l,k}\}$.
The ATC algorithm \eqref{Equ:DiffusionAdaptation:ATC0} involves two steps. In the first step, node $k$ receives
gradient vector information from its neighbors and uses it to update its estimate $w_{k,i-1}$ to an
intermediate value $\psi_{k,i}$. All other nodes in the network are performing a similar step and generating
their intermediate estimate $\psi_{l,i}$. In the second step, node $k$ aggregates the estimates
$\{\psi_{l,i}\}$ of its neighbors and generates $w_{k,i}$. Again, all other nodes are performing
a similar step.
Similarly, if we reverse the order of steps
\eqref{Equ:DiffusionAdaptation:Adaptation_intermedieate}
and \eqref{Equ:DiffusionAdaptation:Combination_intermediate} to implement
\eqref{Equ:DiffusionAdaptation:GradientDescent}, we can motivate the following alternative
Combine-then-Adapt (CTA) diffusion strategy
(whose structure is similar to the CTA algorithm originally proposed in \cite{lopesdistributed,lopes2007diffusion,cattivelli2007diffusionRLS,Sayed07,lopes2008diffusion,
cattivelli2008diffusion,Cattivelli10}  for mean-square-error estimation):
	\begin{align}
		\boxed{
				\label{Equ:DiffusionAdaptation:CTA0}
				\begin{array}{l}
					\psi_{k,i\!-\!1}	=	\displaystyle \sum_{l \in \mathcal{N}_k} a_{l,k} w_{l,i-1}		 \\
					w_{k,i}	=	 \displaystyle\psi_{k,i-1}
								- \mu_k \sum_{l \in \mathcal{N}_k} c_{l,k} \nabla_w J_l(\psi_{k,i-1})
				\end{array}
			}
	\end{align}

\noindent
Adaptive diffusion strategies of the above ATC and CTA types were first proposed and extended in
\cite{lopesdistributed,lopes2007diffusion,cattivelli2007diffusionRLS,Sayed07,lopes2008diffusion,
cattivelli2008diffusion,cattivelli2008TSPdiffusionRLS,Cattivelli10,cattivelli2010TACdiffusionKalman}
for the solution of distributed mean-square-error, least-squares, and state-space estimation
problems over networks.
{\color{black}
The special form of ATC strategy \eqref{Equ:DiffusionAdaptation:ATC0} for minimum-mean-square-error
estimation is listed further ahead as Eq. \eqref{Equ:DiffusionAdaptation:ATC_MMSE} in Example \ref{Example:ATC_MMSE};
the same strategy as \eqref{Equ:DiffusionAdaptation:ATC_MMSE}
also appeared in \cite{stankovic2011decentralized}
albeit with a vanishing step-size sequence to ensure convergence towards consensus.
}%
A special case of the diffusion strategy \eqref{Equ:DiffusionAdaptation:CTA0}
(corresponding to choosing $c_{l,k}=0$ for $l\neq k$ and $c_{k,k}=1$, i.e., without sharing gradient information)
was used in the works \cite{ram2010distributed,bianchi2011convergence,srivastava2011distributed}
to solve distributed
optimization problems that require all nodes to reach agreement
about $w^o$ by relying on step-sizes that decay to zero with time.
Diffusion recursions
of the forms \eqref{Equ:DiffusionAdaptation:ATC0} and \eqref{Equ:DiffusionAdaptation:CTA0}
are more general than these earlier investigations in a couple of respects.
First, they do not only
diffuse the local estimates, but they can also diffuse the local gradient vectors.
In other words, two sets of combination coefficients $\{a_{l,k},c_{l,k}\}$ are used.
Second, the combination weights $\{a_{l,k}\}$ are not required to be doubly stochastic
(which would require both the rows and columns of the weighting matrix $A=[a_{l,k}]$
to add up to one; as seen from \eqref{Equ:AdaptiveDiffusion:Condition_a},
we only require the entries on the columns of $A$ to add up to one).
Finally, and most importantly, the step-size parameters $\{\mu_k\}$ in
\eqref{Equ:DiffusionAdaptation:ATC0} and \eqref{Equ:DiffusionAdaptation:CTA0} are not required to depend
on the time index $i$ and are not required to vanish as $i\rightarrow\infty$. Instead, they can
assume constant values, which is critical to endow the network with \emph{continuous}
adaptation and learning abilities (otherwise, when step-sizes die out, the network stops learning).
Constant step-sizes also endow networks with tracking abilities, in which case the algorithms can track time changes in the optimal $w^o$.

Constant step-sizes will be shown further ahead to be  sufficient to guarantee agreement among the nodes when
there is no noise in the data. However, when measurement noise and gradient noise are present,
using constant step-sizes does not \emph{force} the nodes to attain agreement about $w^o$ (i.e., to converge to the same $w^o$).
 Instead, the nodes will
be shown to
tend to individual estimates for $w^o$ that are within a small mean-square-error (MSE) bound
from the optimal solution; the bound will be proportional to the step-size so that sufficiently small step-sizes lead to small MSE values. Multi-agent
systems in nature behave in this manner; they do not require exact agreement among their agents but allow for
fluctuations due to individual noise levels (see \cite{cattivelli2011modeling, tu2011fish}).
Giving individual nodes this flexibility, rather than forcing them to operate in agreement with the remaining nodes, ends up leading to nodes with enhanced learning abilities.

{
Before proceeding to a detailed analysis of the performance of the diffusion algorithms \eqref{Equ:DiffusionAdaptation:ATC0}--\eqref{Equ:DiffusionAdaptation:CTA0},
we note that these strategies differ in important ways
from traditional consensus-based distributed solutions, which are of the following form
\cite{nedic2009distributed,bertsekasparallel,nedic2009bookchapter,kar2011converegence}:
    \begin{align}
        \label{Equ:DiffusionAdaptation:ConsensusTypeAlgorithm}
        w_{k,i}	=	 \sum_{l \in \mathcal{N}_k} a_{l,k} w_{k,i-1}
								- \mu_k(i) \cdot \nabla_w J_l(w_{k,i-1})
    \end{align}
usually with a time-variant step-size sequence, $\mu_k(i)$, that decays to zero.
For example, if we set $C \triangleq [c_{l,k}] = I$ in the CTA algorithm
\eqref{Equ:DiffusionAdaptation:CTA0}
and substitute the combination step into the adaptation step, we obtain:
    \begin{align}
        \label{Equ:DiffusionAdaptation:CTA_oneline}
        w_{k,i}	=	 \sum_{l \in \mathcal{N}_k} a_{l,k} w_{k,i-1}
								- \mu_k \nabla_w J_l\Big(\sum_{l \in \mathcal{N}_k} a_{l,k} w_{k,i-1}\Big)
    \end{align}
Thus, note that
the gradient vector in \eqref{Equ:DiffusionAdaptation:CTA_oneline}
is evaluated at $\psi_{k,i-1}$, while in
\eqref{Equ:DiffusionAdaptation:ConsensusTypeAlgorithm} it is evaluated at $w_{k,i-1}$.
Since $\psi_{k,i-1}$ already incorporates information from neighbors, we would expect
the diffusion algorithm to perform better.
Actually, it is shown in \cite{Tu2012diffcons} that,
for mean-square-error estimation problems, diffusion strategies achieve
higher convergence rate and lower mean-square-error than consensus strategies
due to these differences in the dynamics of the algorithms.
}

\section{Mean-Square Performance Analysis}
\label{Sec:ConvergenceAnalysis}
The diffusion algorithms \eqref{Equ:DiffusionAdaptation:ATC0} and
\eqref{Equ:DiffusionAdaptation:CTA0}
depend on sharing local gradient vectors $\nabla_{w} J_l(\cdot)$.
In many cases of practical relevance, the exact gradient vectors are not available
and approximations are instead used.
We model the inaccuracy in the gradient vectors as some {\em random}
additive noise component, say, of the form:
	\begin{align}
		\label{Equ:Diffusion:NoisyGradient}
		\widehat{\nabla}_w J_l({w}) = \nabla_w J_l({w}) + \bm{v}_{l,i}({w})
	\end{align}
where $\bm{v}_{l,i}(\cdot)$ denotes the perturbation and is often referred
to as gradient noise. Note that we are using a boldface symbol $\bm{v}$
to refer to the gradient noise  since it is generally stochastic in nature.
\begin{example}
\label{FN:GradientNoiseExample}
{\rm
Assume the individual
cost $J_l(w)$ at node $l$ can be expressed as the expected value of a certain loss function
$Q_l(\cdot,\cdot)$,
i.e., $J_l(w)=\E\{Q_l(w,\bm{x}_{l,i})\}$,
where the expectation is with respect to the randomness in the data
samples $\{\bm{x}_{l,i}\}$ that are collected at node $l$ at time $i$.
Then, if we replace
the true gradient $\nabla_w J_l(w)$ with its stochastic gradient approximation
$\widehat{\nabla}_wJ_l(w)=\nabla_w Q_l(w,\bm{x}_{l,i})$, we find that the gradient noise in this case can
be expressed as
    \begin{align}
        \label{Equ:PerformanceAnalysis:GradientNoise_ExpectedLoss}
        \bm{v}_{l,i}(w)   =   \nabla_w Q_l(w,\bm{x}_{l,i})
                            -
                            \nabla_w \E\{Q_l(w,\bm{x}_{l,i})\}
    \end{align}
}
\hfill\IEEEQED
\end{example}
\noindent
Using the perturbed gradient vectors \eqref{Equ:Diffusion:NoisyGradient},
the diffusion algorithms \eqref{Equ:DiffusionAdaptation:ATC0}--\eqref{Equ:DiffusionAdaptation:CTA0}
become the following:
	\begin{align}
		(\mathrm{ATC})~
		\boxed{
				\label{Equ:DiffusionAdaptation:ATC}
				\begin{array}{l}
					\bm{\psi}_{k,i}	=	\displaystyle
									\bm{w}_{k,i-1}
									\!-\!
									\mu_k \sum_{l \in \mathcal{N}_k}  c_{l,k}
									\widehat{\nabla}_w J_l(\bm{w}_{k,i-1})	\\
					\bm{w}_{k,i}	=	\displaystyle \sum_{l \in \mathcal{N}_k} a_{l,k} \bm{\psi}_{l,i}
				\end{array}
			}
	\end{align}
	\begin{align}
		(\mathrm{CTA})~
		\boxed{
				\label{Equ:DiffusionAdaptation:CTA}
				\begin{array}{l}
					\bm{\psi}_{k,i-1}	=	\displaystyle \sum_{l \in \mathcal{N}_k} a_{l,k} \bm{w}_{l,i-1}	 \\
					\bm{w}_{k,i}	=	\displaystyle \bm{\psi}_{k,i-1}
									\!-\!
									\mu_k \sum_{l \in \mathcal{N}_k}  c_{l,k}
									\widehat{\nabla}_w J_l(\bm{\psi}_{k,i-1})
				\end{array}
			}
	\end{align}

\noindent Observe that, starting with
\eqref{Equ:DiffusionAdaptation:ATC}--\eqref{Equ:DiffusionAdaptation:CTA}, we will be using boldface letters to refer to the various estimate quantities in order to highlight the fact that
 they are also stochastic in nature due to the presence of the gradient noise.

Given the above algorithms, it is necessary to examine their performance in light of the
approximation steps \eqref{Equ:ProblemFormulation:J_glob_approx}--%
\eqref{Equ:DiffusionAdaptation:Combination_intermediate2}
that were employed to arrive at them,
and in light of the gradient noise \eqref{Equ:Diffusion:NoisyGradient} that seeps into the recursions.
A convenient framework
to carry out this analysis is mean-square analysis. In this framework, we assess
how close the individual estimates $\bm{w}_{k,i}$ get to the minimizer $w^o$ in the mean-square-error
(MSE) sense.
{\color{black}
In practice, it is not necessary to force the individual agents to reach agreement
and to converge to the same $w^o$ using diminishing step-sizes. It is sufficient for the nodes to
converge within acceptable MSE bounds from $w^o$. This flexibility is beneficial and
is common in biological networks; it allows nodes to learn and adapt in time-varying environments without the forced requirement of having to agree with neighbors.
}

The main results that we derive in this section are summarized as follows. First, we derive conditions
on the \emph{constant step-sizes} to ensure boundedness and convergence of
the mean-square-error for sufficiently small step-sizes ---
see
\eqref{Equ:ConvergenceAnalysis:MS_Theorem:Stepsizes} and
\eqref{Equ:ConvergenceAnalysis:MSD_k} further ahead.
Second, despite the fact that nodes influence each other's behavior,
we are able to quantify the performance of every node in the network and to
derive
closed-form expressions for the mean-square performance at small step-sizes ---
see \eqref{Equ:ConvergenceAnalysis:MSD_k}--\eqref{Equ:ConvergenceAnalysis:W_infty}.
Finally, as a special case, we are able to show that constant step-sizes can still ensure that the estimates
across all nodes converge to the optimal $w^o$ and reach agreement in the \emph{absence of noise}
--- see Theorem \ref{Corollary:ConvergenceAnalysis:Convergence_NoiseFreeCase}.

Motivated by \cite{Cattivelli10}, we address the mean-square-error performance of the adaptive ATC and CTA
diffusion strategies \eqref{Equ:DiffusionAdaptation:ATC}--\eqref{Equ:DiffusionAdaptation:CTA}
by treating them as special cases of a general
diffusion structure of the following form:
	\begin{align}
		\label{Equ:ConvergenceAnalysis:Diffusion_General}
				\bm{\phi}_{k,i-1}	&=	\displaystyle \sum_{l=1}^N p_{1,l,k} \bm{w}_{l,i-1}		\\
		\label{Equ:ConvergenceAnalysis:Diffusion_General1}
				\bm{\psi}_{k,i}		&=
									\displaystyle
									\bm{\phi}_{k,i-1}
									-
									\mu_k \sum_{l=1}^N
									s_{l,k}
									\big[
										\nabla_w J_l(\bm{\phi}_{k,i-1})
										+
										\bm{v}_{l,i}(\bm{\phi}_{k,i-1})
									\big]
									\\
		\label{Equ:ConvergenceAnalysis:Diffusion_General2}
				\bm{w}_{k,i}		&=	\displaystyle \sum_{l=1}^N p_{2,l,k} \bm{\psi}_{l,i}
	\end{align}
The coefficients $\{p_{1,l,k}\}$, $\{s_{l,k}\}$, and $\{p_{2,l,k}\}$ are nonnegative real
coefficients corresponding to the $\{l,k\}$-th entries of three matrices $P_1$, $S$, and $P_2$, respectively.
Different choices for  $\{P_1,P_2,S\}$ correspond to different cooperation modes. For example,
the choice $P_1 = I$, $P_2 = I$ and $S=I$ corresponds to the non-cooperative case where nodes
do not interact.
On the other hand, the choice $P_1 = I$, $P_2=A=[a_{l,k}]$ and $S=C=[c_{l,k}]$ corresponds to ATC
\cite{lopes2008diffusion,cattivelli2008diffusion,Cattivelli10},
while the choice $P_1 = A$, $P_2 = I$ and $S=C$ corresponds to CTA
\cite{lopesdistributed,lopes2007diffusion,Sayed07,lopes2008diffusion,
cattivelli2008diffusion,Cattivelli10}.  We can
also set $S=I$ in ATC and CTA to derive simplified versions that have no gradient exchange\cite{lopes2008diffusion}.
Furthermore, if in CTA ($P_2=I$), we enforce $P_1=A$ to be doubly stochastic, set $S=I$,
and use a time-decaying step-size parameter ($\mu_k(i)\rightarrow 0$), then
we obtain the unconstrained
version used by \cite{ram2010distributed,srivastava2011distributed}. The matrices
$\{P_1,P_2,S\}$ are required to satisfy:
	\begin{align}
		\label{Equ:ConvergenceAnalysis:ConditionCombinationWeights}
		\boxed{
			P_1^T \mathds{1} = \mathds{1}, \;
			P_2^T \mathds{1} = \mathds{1}, \;
			S \mathds{1} = \mathds{1}
		}
	\end{align}
where the notation $\mathds{1}$ denotes a vector whose entries are all equal to one.

\subsection{Error Recursions}

We first derive the error recursions corresponding to the general diffusion formulation in
\eqref{Equ:ConvergenceAnalysis:Diffusion_General}--\eqref{Equ:ConvergenceAnalysis:Diffusion_General2}.
Introduce the error vectors:
	\begin{align}		
		\label{Equ:ConvergenceAnalysis:Def_errorQuantities}
		\tilde{\bm{\phi}}_{k,i}	\triangleq	w^o\!-\!\bm{\phi}_{k,i},
		\;
		\tilde{\bm{\psi}}_{k,i}	\triangleq	w^o\!-\!\bm{\psi}_{k,i},
		\;
		\tilde{\bm{w}}_{k,i}	\triangleq	w^o\!-\!\bm{w}_{k,i}
	\end{align}
Then, subtracting both sides of
\eqref{Equ:ConvergenceAnalysis:Diffusion_General}--\eqref{Equ:ConvergenceAnalysis:Diffusion_General2}
from $w^o$ gives:
	\begin{align}
		\label{Equ:ConvergenceAnalysis:phi_ki}
				\tilde{\bm{\phi}}_{k,i-1}	&=	\displaystyle\sum_{l=1}^N p_{1,l,k} \tilde{\bm{w}}_{l,i-1}	 
										\\
		\label{Equ:ConvergenceAnalysis:psi_ki}	
				\tilde{\bm{\psi}}_{k,i}		&=	\displaystyle
										\tilde{\bm{\phi}}_{k,i-1}
										+
										\mu_k  \sum_{l=1}^N  s_{l,k}
										\big[
											\nabla_w  J_l(\bm{\phi}_{k,i-1})
											+
											\bm{v}_{l,i}(\bm{\phi}_{k,i-1})
										\big]
										\\
		\label{Equ:ConvergenceAnalysis:w_ki}
				\tilde{\bm{w}}_{k,i}		&=	\displaystyle\sum_{l=1}^N p_{2,l,k} \tilde{\bm{\psi}}_{l,i}	 
	\end{align}
Expression \eqref{Equ:ConvergenceAnalysis:psi_ki} still includes terms that depend on
$\bm{\phi}_{k,i-1}$ and not on the error quantity, $\tilde{\bm{\phi}}_{k,i-1}$.
We can find a relation in terms of
$\tilde{\bm{\phi}}_{k,i-1}$ by calling upon the following result from \cite[p.24]{poliak1987introduction}
for any twice-differentiable function  $f(\cdot)$:
	\begin{align}
		\label{Equ:ConvergenceAnalysis:MeanValueTheorem_VecVariant}
		\nabla f(y)	=	\nabla f(x) + \left[\int_{0}^{1} \nabla^2 f\big(x\!+\!t(y\!-\!x)\big)dt\right] (y-x)
	\end{align}
where $\nabla^2 f(\cdot)$ denotes the Hessian matrix of $f(\cdot)$ and
is symmetric.
Now, since each component function $J_l(w)$ has a minimizer at $w^o$,
then, $\nabla_w J_l(w^o)=	0$ for $l=1,2,\ldots,N$. Applying
\eqref{Equ:ConvergenceAnalysis:MeanValueTheorem_VecVariant}
to $J_l(w)$ using $x=w^o$ and $y=\bm{\phi}_{k,i-1}$, we get
	\begin{align}
		&\nabla_w J_l(\bm{\phi}_{k,i-1})	
                                            \nonumber\\
        &\quad= \;\nabla_w J_l(w^o)
									-
									\left[
										\int_{0}^{1}
										\nabla_w^2
										J_l\big(w^o - t\tilde{\bm{\phi}}_{k,i-1}\big) dt
									\right]
									\tilde{\bm{\phi}}_{k,i-1}
									\nonumber\\
		\label{Equ:ConvergenceAnalysis:Gradient_LinearRepresentation}
								&\quad\triangleq	-
									\bm{H}_{l,k,i-1}
									\tilde{\bm{\phi}}_{k,i-1}
	\end{align}
where we are introducing the symmetric random matrix
	\begin{align}
		\label{Equ:ConvergenceAnalysis:H_kim1}
		\boxed{
		\bm{H}_{l,k,i-1}	\triangleq	\int_{0}^{1}
							\nabla_w^2
							J_l\big(w^o - t\tilde{\bm{\phi}}_{k,i-1}\big) dt
		}
	\end{align}
Observe that one such matrix is associated with every edge linking two nodes $(l,k)$;
observe further that this matrix changes with time since it depends on the estimate at node $k$.
Substituting
\eqref{Equ:ConvergenceAnalysis:Gradient_LinearRepresentation}--\eqref{Equ:ConvergenceAnalysis:H_kim1}
into \eqref{Equ:ConvergenceAnalysis:psi_ki} leads to:
	\begin{align}		
				\tilde{\bm{\psi}}_{k,i}		=\;&	\displaystyle
											\Big[
												I_M - \mu_k \sum_{l=1}^N s_{l,k} \bm{H}_{l,k,i-1}
											\Big]
											\tilde{\bm{\phi}}_{k,i-1}
                                            \nonumber\\
		\label{Equ:ConvergenceAnalysis:ErrorRecursions}	
											\;&+
											\mu_k
											\sum_{l=1}^N
											s_{l,k}
											\bm{v}_{l,i}(\bm{\phi}_{k,i-1})
	\end{align}

We introduce
the network error vectors, which collect the error quantities across all nodes:
	\begin{align}
		\label{Equ:ConvergenceAnalysis:global_error_vector}
		\tilde{\bm{\phi}}_{i}	\triangleq	
							\begin{bmatrix}
								\tilde{\bm\phi}_{1,i}	\\
								\vdots			\\
								\tilde{\bm\phi}_{N,i}
							\end{bmatrix},	\qquad
		\tilde{\bm\psi}_{i}	\triangleq	
							\begin{bmatrix}
								\tilde{\bm\psi}_{1,i}	\\
								\vdots			\\
								\tilde{\bm\psi}_{N,i}
							\end{bmatrix},	\qquad
		\tilde{\bm w}_{i}		\triangleq	
							\begin{bmatrix}
								\tilde{\bm w}_{1,i}	\\
								\vdots			\\
								\tilde{\bm w}_{N,i}
							\end{bmatrix}
	\end{align}
and the following block matrices:
	{
	\begin{align}
		\label{Equ:ConvergenceAnalysis:P1_P2}
		\mathcal{P}_1	\;=\;&	P_1 \otimes I_M,\;
		\mathcal{P}_2	\;=\;	P_2 \otimes I_M	\\
		\label{Equ:ConvergenceAnalysis:S_M}
		\mathcal{S}	\;=\;&	S \otimes I_M, \;
		\mathcal{M}	\;=\;	\Omega \otimes I_M \\
		\label{Equ:ConvergenceAnalysis:Xi}
		\Omega			\;=\;&	\mathrm{diag}
						\left\{
							\mu_1, \; \ldots, \; \mu_N
						\right\}
						\\
		\label{Equ:ConvergenceAnalysis:D_i_minus_1}
		\!\!\!\bm{\mathcal{D}}_{i-1}	
					\;=\;&	\sum_{l=1}^N
						\mathrm{diag}
						\big\{
							 s_{l,1} \bm{H}_{l,1,i-1},
							\cdots,
							s_{l,N} \bm{H}_{l,N,i-1}
						\big\}
						\\
		\label{Equ:ConvergenceAnalysis:G_i}
		\bm{g}_i
					\;=\;&	\sum_{l=1}^N\!
						\mathrm{col}
						\big\{
							 s_{l,1}\bm{v}_{l,i}(\bm{\phi}_{1,i\!-\!1}),
							\cdots,\!
							s_{l,N}\bm{v}_{l,i}(\bm{\phi}_{N,i\!-\!1})
						\big\}
	\end{align}
	}%
where the symbol $\otimes$ denotes Kronecker products \cite{horn1990matrix}. Then, recursions
\eqref{Equ:ConvergenceAnalysis:phi_ki}, \eqref{Equ:ConvergenceAnalysis:ErrorRecursions}
and \eqref{Equ:ConvergenceAnalysis:w_ki} lead to:
	\begin{align}
		\label{Equ:ConvergenceAnalysis:ErrorRecursion_final}
		\boxed{
			\tilde{\bm w}_i	=	\mathcal{P}_2^T
							[I_{MN}-\mathcal{M}\bm{\mathcal{D}}_{i-1}]
							\mathcal{P}_1^T
							\tilde{\bm w}_{i-1}
							+
							\mathcal{P}_2^T \mathcal{M} \bm{g}_i
		}
	\end{align}
To proceed with the analysis, we introduce the following assumption on the cost functions and gradient noise,
followed by a lemma on $\bm{H}_{l,k,i-1}$.
	\begin{assumption}[Bounded Hessian]
		\label{Assumption:Hessian}
		Each component cost function $J_l(w)$ has a bounded Hessian matrix, i.e.,
		there exist
		nonnegative real numbers $\lambda_{l,\min}$ and $\lambda_{l,\max}$
        such that $\lambda_{l,\min} \le \lambda_{l,\max}$ and that
        for all $w$:
			\begin{align}
				\label{Assumption:StrongConvexity}
				\lambda_{l,\min} I_M	\le		\nabla_w^2 J_l(w)		\le	\lambda_{l,\max} I_M
			\end{align}
		Furthermore, the $\{\lambda_{l,\min}\}_{l=1}^N$ satisfy
			\begin{align}
				\label{Assumption:StrongConvexity1}
				\sum_{l=1}^N s_{l,k} \lambda_{l,\min} >0, k=1,2,\ldots,N
			\end{align}
				\hfill\IEEEQED
	\end{assumption}
\noindent Condition \eqref{Assumption:StrongConvexity1}
ensures that the local cost functions $\{J_{k}^{\mathrm{loc}} (w)\}$ defined earlier in
\eqref{Equ:ProblemFormulation:J_k_loc0} are strongly convex and, hence,
have a unique minimizer at $w^o$.
	\begin{assumption}[Gradient noise]
		\label{Assumption:GradientNoise}
        There exist $\alpha \ge 0$ and $\sigma_{v}^2 \ge 0$ such that,
        for all $\bm{w} \in \mc{F}_{i-1}$ and for all $i$, $l$:
		{
			\begin{align}
				\label{Assumption:GradientNoise:ZeroMean_Uncorrelated}
				&\mathbb{E}\left\{\bm{v}_{l,i}(\bm{w}) \;|\; \mc{F}_{i-1} \right\} = 0	 \\
				\label{Assumption:GradientNoise:Norm}
				&\mathbb{E}\left\{ \|\bm{v}_{l,i}(\bm{w}) \|^2 \right\}
						\le \alpha \E\|w^o-\bm{w}\|^2 + \sigma_{v}^2
			\end{align}
		}%
        where  $\mc{F}_{i-1}$  denotes  the  past  history ($\sigma-$field)  of  estimates
        $\{\bm{w}_{k,j}\}$ for $j \le i-1$ and all $k$.
		\hfill\IEEEQED
	\end{assumption}	
	\begin{lemma}[Bound on $\bm{H}_{l,k,i-1}$]
		\label{Lemma:H_lki_bounds}
		Under Assumption \ref{Assumption:Hessian},
		the matrix $\bm{H}_{l,k,i-1}$ defined in \eqref{Equ:ConvergenceAnalysis:H_kim1}
		is a nonnegative-definite matrix that satisfies:
			\begin{align}
				\label{Equ:ConvergenceAnalysis:H_lki_Bounds}
				\lambda_{l,\min} I_M	\le	\bm{H}_{l,k,i-1}	\le	\lambda_{l,\max} I_M
			\end{align}
	\end{lemma}
	\begin{IEEEproof}
		It suffices to prove that $\lambda_{l,\min} \le x^T \bm{H}_{l,k,i-1} x \le \lambda_{l,\max}$
		for arbitrary $M \times 1$ unit Euclidean norm vectors $x$.
		By \eqref{Equ:ConvergenceAnalysis:H_kim1}  and \eqref{Assumption:StrongConvexity},  we have
			\begin{align}
				x^T \bm{H}_{l,k,i-1} x	
							&=	\int_{0}^{1}
                                x^T\nabla_w^2 J_l\left(w^o - t\tilde{\bm{\phi}}_{k,i-1}\right) x \; dt
                                \nonumber\\
							&\le	\int_0^1 \lambda_{l,\max} dt = \lambda_{l,\max}
								\nonumber
			\end{align}
		In a similar way, we can verify that $x^T \bm{H}_{l,k,i-1} x \ge \lambda_{l,\min}$.
	\end{IEEEproof}	
    \vspace{0.5em}

        {\color{black}
		In distributed subgradient methods (e.g.,
        \cite{nedic2009distributed,ram2010distributed,srivastava2011distributed}), the norms of the subgradients
        are usually required to be uniformly bounded. Such assumption is restrictive in the
        unconstrained optimization of differentiable functions.
        }
		Assumption \ref{Assumption:Hessian} is more relaxed in that
		it allows the gradient vector $\nabla_w J_l(w)$ to have unbounded norm
        (e.g., quadratic costs).
		Furthermore, condition \eqref{Assumption:GradientNoise:Norm}
		allows the variance of the gradient noise
		to grow no faster than
		$\E\|w^o-\bm{w}\|^2$.
		This condition is also more general than the
		uniform bounded assumption used in \cite{ram2010distributed}
		(Assumptions 5.1 and 6.1), which requires instead:
			\begin{align}
                \label{Equ:AbsoluteRandomNoise_GradientNoise}
				\!\!\!\mathbb{E} \| \bm{v}_{l,i}(\bm{w}) \|^2 \!\le\!	\sigma_v^2,
				\quad
				\mathbb{E}
                \left\{ \|\bm{v}_{l,i}(\bm{w}) \|^2  |  \mathcal{F}_{i\!-\!1} \right\}
                \!\le\!	\sigma_v^2
			\end{align}
		Furthermore, condition \eqref{Assumption:GradientNoise:Norm} is similar to
		condition (4.3) in \cite[p.635]{bertsekas2000gradient}:
			\begin{align}
				\label{Equ:ConvergenceAnalysis:GradientNoise_Bertsekas}
				\mathbb{E}\left\{ \| \bm{v}_{l,i}(\bm{w}) \|^2 | \mathcal{F}_{i-1} \right\}
						\le \alpha \big[\|\nabla_w J_l(\bm{w})\|^2 + 1\big]
			\end{align}
		which is a combination of the ``relative random noise''  and the
		``absolute random noise'' conditions defined in \cite[pp.100--102]{poliak1987introduction}.
		Indeed, we can derive \eqref{Assumption:GradientNoise:Norm}
		by substituting \eqref{Equ:ConvergenceAnalysis:Gradient_LinearRepresentation}
		into \eqref{Equ:ConvergenceAnalysis:GradientNoise_Bertsekas},
        taking expectation with respect to $\mathcal{F}_{i-1}$,
		and then using \eqref{Equ:ConvergenceAnalysis:H_lki_Bounds}.
{
\color{black}
\begin{example}
    \rm
    Such a mix of ``relative random noise'' and ``absolute random noise'' is of practical importance.
    For instance, consider an example in which the loss function at node $l$ is chosen to be
    of the following quadratic form:
        \begin{align}
            Q_l(w,\{\bm{u}_{l,i},\bm{d}_l(i)\})=|\bm{d}_l(i)-\bm{u}_{l,i}w|^2
                                                    \nonumber
        \end{align}
    for some scalars $\{\bm{d}_l(i)\}$ and $1\times M$ regression vectors $\{\bm{u}_{l,i}\}$. The corresponding cost function is then:
        \begin{align}
            \label{Equ:PerformanceAnalysis:QuadraticCost_LMS}
            J_l(w)  &=   \E |\bm{d}_l(i)-\bm{u}_{l,i}w|^2
        \end{align}
    Assume further that the data $\{\bm{u}_{l,i},\bm{d}_l(i)\}$ satisfy the linear regression model
        \begin{align}
            \label{Equ:PerformanceAnalysis:gradientNoise_linearmodel}
            \bm{d}_l(i) =   \bm{u}_{l,i}w^o + \bm{z}_l(i)
        \end{align}
    where the regressors $\{\bm{u}_{l,i}\}$ are zero mean and independent over time with
    covariance matrix $R_{u,l}=\E \{\bm{u}_{l,i}^T\bm{u}_{l,i}\}$, and
    the noise sequence $\{\bm{z}_k(j)\}$ is also zero mean, white, with variance $\sigma_{z,k}^2$,
    and independent of the regressors $\{\bm{u}_{l,i}\}$ for all $l,k,i,j$.
    Then, using
    \eqref{Equ:PerformanceAnalysis:gradientNoise_linearmodel} and
    \eqref{Equ:PerformanceAnalysis:GradientNoise_ExpectedLoss},
    the gradient noise in this case can be expressed as:
        \begin{align}
            \label{Equ:PerformanceAnalysis:GradientNoise_LMS}
            \bm{v}_{l,i}(\bm{w})
                =   2(R_{u,l}-\bm{u}_{l,i}^T\bm{u}_{l,i})(w^o-\bm{w}) - 2\bm{u}_{l,i}^T \bm{z}_l(i)
        \end{align}
    It can easily be verified that this noise satisfies both conditions stated in Assumption
    \ref{Assumption:GradientNoise}, namely, \eqref{Assumption:GradientNoise:ZeroMean_Uncorrelated}
    and also:
        \begin{align}
            &\mathbb{E}\left\{ \|\bm{v}_{l,i}(\bm{w})\|^2\right\}
                                    \nonumber\\
            \label{Equ:PerformanceAnalysis:GradientNoiseNorm_LMS}
    						&\le 4  \E\|R_{u,l}\!-\!\bm{u}_{l,i}^T\bm{u}_{l,i}\|^2
                                    \cdot
                                    \E\|w^o\!-\!\bm{w}\|^2 \!+\! 4 \sigma_{z,l}^2 \mathrm{Tr}(R_{u,l})
        \end{align}
    for all $\bm{w} \in \mc{F}_{i-1}$.
    Note that both relative random noise and absolute random noise components
    appear in \eqref{Equ:PerformanceAnalysis:GradientNoiseNorm_LMS} and
    are necessary to model the statistical gradient perturbation even for quadratic costs.
    Such costs, and linear regression models of the form
    \eqref{Equ:PerformanceAnalysis:gradientNoise_linearmodel},
    arise frequently in the context of adaptive filters --- see, e.g.,
    \cite{Sayed08,haykin2002adaptive,lopes2007incremental,lopesdistributed,lopes2007diffusion,
    Sayed07,lopes2008diffusion,cattivelli2008diffusion,Cattivelli10,cattivelli2007diffusionRLS,cattivelli2008TSPdiffusionRLS,
    cattivelli2011modeling,arenas2006plant,silva2008improving,theodoridis2011adaptive}.
    \hfill\IEEEQED
\end{example}
}
\begin{example}	
\label{Example:ATC_MMSE}
\rm
Quadratic costs of the form \eqref{Equ:PerformanceAnalysis:QuadraticCost_LMS} are common in mean-square-error estimation for linear regression models of the type
\eqref{Equ:PerformanceAnalysis:gradientNoise_linearmodel}. If we use the
instantaneous approximations as is common in the context of stochastic approximation
and adaptive filtering \cite{Sayed08,haykin2002adaptive,poliak1987introduction},
then the actual gradient $\nabla_w J_l(w)$ can be approximated by
    \begin{align}
        \widehat{\nabla}_{w} J_{l}(w)   &=   \nabla_wQ_l(w,\{\bm{u}_{l,i},\bm{d}_l(i)\})      \nonumber\\
                                        &=   -2\bm{u}_{l,i}^T[\bm{d}_{l}(i)-\bm{u}_{l,i} w]
    \end{align}
Substituting into \eqref{Equ:DiffusionAdaptation:ATC}--\eqref{Equ:DiffusionAdaptation:CTA},
and assuming $C=I$ for illustration purposes only,
we arrive at the following ATC and CTA diffusion strategies
originally proposed and extended in \cite{lopesdistributed,lopes2007diffusion,
Sayed07,lopes2008diffusion, cattivelli2008diffusion,Cattivelli10}
for the solution of distributed mean-square-error estimation problems:
    \begin{align}
		\!\!(\mathrm{ATC})~\boxed{
				\label{Equ:DiffusionAdaptation:ATC_MMSE}
				\begin{array}{l}
					\bm{\psi}_{k,i}	=	\displaystyle
									\bm{w}_{k,i\!-\!1}
									\!+\!
									2\mu_k
									\bm{u}_{k,i}^T[\bm{d}_{k}(i)\!-\!\bm{u}_{k,i} \bm{w}_{k,i\!-\!1}]	\\
					\bm{w}_{k,i}	=	\displaystyle \sum_{l \in \mathcal{N}_k} a_{l,k} \bm{\psi}_{l,i}
				\end{array}
			}
	\end{align}
	\begin{align}
		\!\!(\mathrm{CTA})~\boxed{
				\label{Equ:DiffusionAdaptation:CTA_MMSE}
				\begin{array}{l}
					\bm{\psi}_{k,i-1}	=	\displaystyle \sum_{l \in \mathcal{N}_k} a_{l,k} \bm{w}_{l,i-1}	 \\
					\bm{w}_{k,i}	=	\displaystyle \bm{\psi}_{k,i\!-\!1}
									\!+\!
									2\mu_k \bm{u}_{k,i}^T[\bm{d}_{k}(i)
                                    \!-\!
                                    \bm{u}_{k,i} \bm{\psi}_{k,i\!-\!1}]
				\end{array}
			}
	\end{align}	
\hfill\IEEEQED
\end{example}

\subsection{Variance Relations}
\label{Sec:ConvergenceAnalysis:VarianceRelation}
The purpose of the mean-square analysis in the sequel
is to answer two questions in the presence of gradient perturbations. First, how small the mean-square error,
$\mathbb{E}\|\tilde{\bm{w}}_{k,i}\|^2$, gets as $i\rightarrow\infty$ for any of the nodes $k$. Second, how fast this error
variance tends towards its steady-state value. The first question pertains to steady-state performance and the second
question pertains to transient/convergence rate performance. Answering such questions for a distributed algorithm
over a network is a challenging task largely because the nodes influence each other's behavior: performance at one
node diffuses through the network to the other nodes as a result of the topological constraints linking the nodes.
The approach we take to examine the mean-square performance of the diffusion algorithms is by studying how the
variance $\mathbb{E}\|\tilde{\bm{w}}_{k,i}\|^2$, or a weighted version of it, evolves over time. As the derivation will
show, the evolution of this variance satisfies a nonlinear relation. Under some
reasonable assumptions on the noise profile, and the local cost functions, we will
be able to bound these error variances as well as estimate their steady-state values for sufficiently small step-sizes.
We will also derive closed-form expressions that characterize the network performance.
The details are as follows.

Equating the squared \emph{weighted} Euclidean norm of both sides of
\eqref{Equ:ConvergenceAnalysis:ErrorRecursion_final}, applying the expectation operator
and using using \eqref{Assumption:GradientNoise:ZeroMean_Uncorrelated},
we can show that the following variance relation holds:
	\begin{equation}
		\label{Equ:ConvergenceAnalysis:WeightedEnergyConservation_Relation}
		\boxed{
			\begin{split}
			&\mathbb{E}\|\tilde{\bm{w}}_{i}\|_\Sigma^2
						=	\mathbb{E}\big\{\|\tilde{\bm{w}}_{i-1}\|_{\bm{\Sigma}'}^2\big\}
							+
							\mathbb{E}\|\mathcal{P}_2^T \mathcal{M} \bm{g}_i\|_\Sigma^2
							\\
			&\bm{\Sigma}'	=	\mathcal{P}_1
							[I_{MN}\!-\!\mathcal{M}\bm{\mathcal{D}}_{i\!-\!1}]
							\mathcal{P}_2
							\Sigma
							\mathcal{P}_2^T
							[I_{MN}\!-\!\mathcal{M}\bm{\mathcal{D}}_{i\!-\!1}]
							\mathcal{P}_1^T
			\end{split}
		}
	\end{equation}
where $\Sigma$ is a positive semi-definite weighting matrix that we are free to choose.
The variance expression \eqref{Equ:ConvergenceAnalysis:WeightedEnergyConservation_Relation}
shows how the quantity $\mathbb{E}\|\tilde{\bm{w}}_i\|^2_{\Sigma}$ evolves with time. Observe, however, that the weighting matrix on $\tilde{\bm{w}}_{i-1}$ on the right-hand side of
\eqref{Equ:ConvergenceAnalysis:WeightedEnergyConservation_Relation}
is a different matrix, denoted by $\bm{\Sigma}'$, and this matrix is actually random in nature
(while $\Sigma$ is deterministic).
As such, result \eqref{Equ:ConvergenceAnalysis:WeightedEnergyConservation_Relation}
is not truly a recursion. Nevertheless, it is possible, under a small step-size approximation, to rework
variance relations such as
\eqref{Equ:ConvergenceAnalysis:WeightedEnergyConservation_Relation}
into a recursion by following certain steps that are characteristic of the energy conservation approach to
mean-square analysis \cite{Sayed08}.

The first step in this regard would be to replace $\bm{\Sigma}'$ by its mean $\mathbb{E}\bm{\Sigma}'$.
However, the matrix $\bm{\Sigma}'$ depends
on the $\{\bm{H}_{l,k,i-1}\}$ via $\bm{\mathcal{D}}_{i-1}$
(see \eqref{Equ:ConvergenceAnalysis:D_i_minus_1}).
It follows from the definition of $\bm{H}_{l,k,i-1}$ in \eqref{Equ:ConvergenceAnalysis:H_kim1}
that $\bm{\Sigma}'$
is dependent on $\tilde{\bm{\phi}}_{k,i-1}$ as well, which in turn is a linear combination of the $\{\tilde{\bm{w}}_{l,i-1}\}$.
Therefore, the main challenge to continue from
\eqref{Equ:ConvergenceAnalysis:WeightedEnergyConservation_Relation}
is that $\bm{\Sigma}'$ depends on $\tilde{\bm{w}}_{i-1}$.
For this reason, we cannot apply directly the traditional step of replacing $\bm{\Sigma}'$ in the first equation of
\eqref{Equ:ConvergenceAnalysis:WeightedEnergyConservation_Relation}
by $\mathbb{E}\bm{\Sigma}'$ as
is typically done in the study of stand-alone adaptive filters to analyze their transient behavior
\cite[p.345]{Sayed08}; in the case of conventional adaptive filters,
the matrix $\bm{\Sigma}'$ is independent of $\tilde{\bm{w}}_{i-1}$.
To address this difficulty, we shall adjust the argument to rely on a set of \emph{inequality}
recursions that will enable us
to bound the steady-state mean-square-error
at each node --- see Theorem \ref{Thm:ConvergenceAnalysis:MS_Theorem}
further ahead.

The procedure is as follows.
First, we note that
$\| x \|^2$ is a convex function of $x$, and that expressions \eqref{Equ:ConvergenceAnalysis:phi_ki}
and \eqref{Equ:ConvergenceAnalysis:w_ki} are convex combinations
of $\{\tilde{\bm{w}}_{l,i-1}\}$ and $\{\tilde{\bm{\psi}}_{l,i}\}$, respectively.
Then, by Jensen's inequality\cite[p.77]{boyd2004convex} and taking expectations, we obtain
{
	\begin{align}
		\label{Equ:ConvergenceAnalysis:Combination1_bound}
		\mathbb{E}\|\tilde{\bm{\phi}}_{k,i-1}\|^2
									&\le		\sum_{l=1}^N p_{1,l,k}\mathbb{E}\| \tilde{\bm{w}}_{l,i-1}\|^2
											\\
		\label{Equ:ConvergenceAnalysis:Combination2_bound}
		\mathbb{E}\|\tilde{\bm{w}}_{k,i}\|^2	&\le		\sum_{l=1}^N p_{2,l,k}\mathbb{E}\| \tilde{\bm{\psi}}_{l,i}\|^2
	\end{align}
}%
for $k=1,\ldots,N$.
Next, we derive a variance relation for \eqref{Equ:ConvergenceAnalysis:ErrorRecursions}.
Equating the squared Euclidean norms of both sides of
\eqref{Equ:ConvergenceAnalysis:ErrorRecursions}, applying the expectation operator, and using
\eqref{Assumption:GradientNoise:ZeroMean_Uncorrelated} from Assumption \ref{Assumption:GradientNoise},
we get
	\begin{align}
		\label{Equ:ConvergenceAnalysis:VarianceRelation_Original}
		\!\!\mathbb{E}\|\tilde{\bm{\psi}}_{k,i}\|^2	&\!= \mathbb{E}\big\{\|\tilde{\bm{\phi}}_{k,i-1}\|_{\bm{\Sigma}_{k,i-1}}^2\big\}
											\!+\!
											\mu_k^2
											\mathbb{E}
											\Big\|
												\sum_{l=1}^N\!
												s_{l,k}
												\bm{v}_{l,i}(\bm{\phi}_{k,i-1})
											\Big\|^2
	\end{align}
where
	\begin{align}
		\bm{\Sigma}_{k,i-1}			
		\label{Equ:ConvergenceAnalysis:VarianceRelation_Sigma_kim1}
								=&\;																										 \Big[
												I_M \!-\! \mu_k \sum_{l=1}^N s_{l,k} \bm{H}_{l,k,i-1}
											\Big]^2
	\end{align}	
We call upon the following two lemmas to bound
\eqref{Equ:ConvergenceAnalysis:VarianceRelation_Original}.
	\begin{lemma}[Bound on $\bm{\Sigma}_{k,i-1}$]
		\label{Lemma:Sigma_k_im1}
		The weighting matrix $\bm{\Sigma}_{k,i-1}$ defined in
		\eqref{Equ:ConvergenceAnalysis:VarianceRelation_Sigma_kim1}
		is a symmetric, positive semi-definite matrix, and satisfies:
			\begin{align}
				\label{Equ:ConvergenceAnalysis:Lemma:Sigma_kim1_Bound}
				0	\le	\bm{\Sigma}_{k,i-1}	\le	\gamma_k^2 I_M
			\end{align}
		where
			\begin{align}
				\label{Equ:ConvergenceAnalysis:gamma_k}
				\boxed{
				\gamma_k		\triangleq		\max
										\Big\{
											\Big|1\!-\!\mu_k \sum_{l=1}^N s_{l,k} \lambda_{l,\max}
											\Big|,\;
											\Big|1\!-\!\mu_k \sum_{l=1}^N s_{l,k} \lambda_{l,\min}\Big|
										\Big\}
				}
			\end{align}
	\end{lemma}
	\begin{IEEEproof}
		By definition \eqref{Equ:ConvergenceAnalysis:VarianceRelation_Sigma_kim1} and the
		fact that $\bm{H}_{l,k,i-1}$ is symmetric --- see
		definition \eqref{Equ:ConvergenceAnalysis:H_kim1}, the matrix
		$I_M \!-\! \mu_k \sum_{l=1}^N s_{l,k} \bm{H}_{l,k,i-1}$ is also symmetric.
		Hence, its square, $\bm{\Sigma}_{k,i-1}$, is symmetric and also nonnegative-definite.
		To establish \eqref{Equ:ConvergenceAnalysis:Lemma:Sigma_kim1_Bound},
		we first use \eqref{Equ:ConvergenceAnalysis:H_lki_Bounds} to note that:
			\begin{align}
				\label{Equ:ConvergenceAnalysis:Proof_Lemma_SimgaBound_Ineq1}
				&I_M \!-\! \mu_k \sum_{l=1}^N s_{l,k} \bm{H}_{l,k,i-1}
					\ge		\Big(1\!-\!\mu_k \sum_{l=1}^N s_{l,k} \lambda_{l,\max}\Big) I_M
							\\
				\label{Equ:ConvergenceAnalysis:Proof_Lemma_SimgaBound_Ineq2}
				&I_M \!-\! \mu_k \sum_{l=1}^N s_{l,k} \bm{H}_{l,k,i-1}
					\le		\Big(1\!-\!\mu_k \sum_{l=1}^N s_{l,k} \lambda_{l,\min}\Big) I_M
			\end{align}
		The matrix $I_M - \mu_k \sum_{l=1}^N s_{l,k} \bm{H}_{l,k,i-1}$ may not be positive semi-definite
		because we have not specified a range for $\mu_k$ yet; the
		expressions on the right-hand side of
		\eqref{Equ:ConvergenceAnalysis:Proof_Lemma_SimgaBound_Ineq1}--%
		\eqref{Equ:ConvergenceAnalysis:Proof_Lemma_SimgaBound_Ineq2} may still be negative.
		However, inequalities
		\eqref{Equ:ConvergenceAnalysis:Proof_Lemma_SimgaBound_Ineq1}--%
		\eqref{Equ:ConvergenceAnalysis:Proof_Lemma_SimgaBound_Ineq2}
		imply that the eigenvalues of $I_M - \mu_k \sum_{l=1}^N s_{l,k} \bm{H}_{l,k,i-1}$
		are bounded as:
			\begin{align}
				\label{Equ:ConvergenceAnalysis:Proof_Lemma_SimgaBound_Ineq3}
                &\lambda\Big(I_M - \mu_k \sum_{l=1}^N s_{l,k} \bm{H}_{l,k,i-1}\Big)
                \ge
				1-\mu_k \sum_{l=1}^N s_{l,k} \lambda_{l,\max}
                \\
                \label{Equ:ConvergenceAnalysis:Proof_Lemma_SimgaBound_Ineq3_2}
				&\lambda\Big(I_M - \mu_k \sum_{l=1}^N s_{l,k} \bm{H}_{l,k,i-1}\Big)
				\le
				1-\mu_k \sum_{l=1}^N s_{l,k} \lambda_{l,\min}
			\end{align}
		By definition \eqref{Equ:ConvergenceAnalysis:VarianceRelation_Sigma_kim1},
		$\bm{\Sigma}_{k,i-1}$ is the square of the symmetric matrix
		$I_M \!-\! \mu_k \sum_{l=1}^N s_{l,k} \bm{H}_{l,k,i-1}$, meaning that
			\begin{align}
				\label{Equ:ConvergenceAnalysis:Proof_Lemma_SimgaBound_eq1}
				\lambda\left(\bm{\Sigma}_{k,i-1}\right)
							=	\left[\lambda\left(I_M - \mu_k \sum_{l=1}^N s_{l,k} \bm{H}_{l,k,i-1}\right)\right]^2
							\ge	0
			\end{align}
		Substituting \eqref{Equ:ConvergenceAnalysis:Proof_Lemma_SimgaBound_Ineq3}--%
        \eqref{Equ:ConvergenceAnalysis:Proof_Lemma_SimgaBound_Ineq3_2}
		into \eqref{Equ:ConvergenceAnalysis:Proof_Lemma_SimgaBound_eq1} leads to
			\begin{align}
				&\lambda\left(\bm{\Sigma}_{k,i-1}\right)	\nonumber\\
				\label{Equ:ConvergenceAnalysis:Proof_Lemma_SimgaBound_Ineq4}
                                            &\le	\max
													\Big\{
														\Big|
															1
															-
															\mu_k \sum_{l=1}^N s_{l,k} \lambda_{l,\max}
														\Big|^2
														,
														\Big|
															1
															-
															\mu_k \sum_{l=1}^N s_{l,k} \lambda_{l,\min}
														\Big|^2
													\Big\}
			\end{align}
		which is equivalent to \eqref{Equ:ConvergenceAnalysis:Lemma:Sigma_kim1_Bound}.
	\end{IEEEproof}
	\begin{lemma}[Bound on noise combination]
		The second term on the right-hand-side of
		\eqref{Equ:ConvergenceAnalysis:VarianceRelation_Original} satisfies:
			\begin{align}					
				\label{Equ:ConvergenceAnalysis:NoiseTerm_bound_final}
					\mathbb{E}
					\Big\|
						\!
						&\sum_{l=1}^N
						s_{l,k}
						\bm{v}_{l,i}(\bm{\phi}_{k,i\!-\!1})
					\Big\|^2		\le		\|S\|_1^2
										\cdot
										\left[
											\alpha \mathbb{E}\|\tilde{\bm{\phi}}_{k,i-1}\|^2
                                            \!+\! \sigma_{v}^2
										\right]
				\end{align}
		where $\|S\|_1$ denotes the $1$-norm of the matrix $S$
		(i.e., the maximum absolute column sum).
	\end{lemma}
	\begin{IEEEproof}
		Applying Jensen's inequality\cite[p.77]{boyd2004convex}, it holds that
			\begin{align}
				&\mathbb{E}
				\Big\|
					\!
					\sum_{l=1}^N
					s_{l,k}
					\bm{v}_{l,i}(\bm{\phi}_{k,i-1})
				\Big\|^2		
                                    \nonumber\\
							&\qquad=		\big(\sum_{l=1}^N s_{l,k}\big)^2 \cdot
									\mathbb{E}
									\Big\|
										\!
										\sum_{l=1}^N
										\frac{s_{l,k}}{\sum_{l=1}^N s_{l,k}}
										\bm{v}_{l,i}(\bm{\phi}_{k,i-1})
									\Big\|^2
									\nonumber\\			
							&\qquad\le		
                                    \big(\sum_{l=1}^N s_{l,k}\big)^2	\cdot
									\sum_{l=1}^N
									\frac{s_{l,k}}{\sum_{l=1}^N s_{l,k}}		
									\mathbb{E}
									\|
										\bm{v}_{l,i}(\bm{\phi}_{k,i-1})
									\|^2
                                    \nonumber\\
							&\qquad=		
                                    \big(\sum_{l=1}^N s_{l,k}\big) \cdot
									\sum_{l=1}^N
									s_{l,k}
									\mathbb{E}
									\|
										\bm{v}_{l,i}(\bm{\phi}_{k,i-1})
									\|^2
                                    \nonumber\\
				\label{Equ:ConvergenceAnalysis:NoiseTerm_bound}
                            &\qquad\le		\big(
										\sum_{l=1}^N
										s_{l,k}
									\big)^2
                                    \cdot
									\big[
									\alpha \mathbb{E}\|\tilde{\bm{\phi}}_{k,i-1}\|^2 + \sigma_{v}^2
									\big]
                                    \\
                \label{Equ:ConvergenceAnalysis:NoiseTerm_bound2}
							&\qquad\le		\|S\|_1^2 \cdot
									\left[\alpha \mathbb{E}\|\tilde{\bm{\phi}}_{k,i-1}\|^2 + \sigma_{v}^2\right]
			\end{align}
        where inequality \eqref{Equ:ConvergenceAnalysis:NoiseTerm_bound} follows by substituting
        \eqref{Assumption:GradientNoise:Norm}, and \eqref{Equ:ConvergenceAnalysis:NoiseTerm_bound2}
        is obtained using the fact that $\|S\|_1$ is the maximum absolute column sum
        and that the entries $\{s_{l,k}\}$ are nonnegative.
	\end{IEEEproof}
    \vspace{1em}
\noindent
Substituting \eqref{Equ:ConvergenceAnalysis:Lemma:Sigma_kim1_Bound} and
\eqref{Equ:ConvergenceAnalysis:NoiseTerm_bound_final} into
\eqref{Equ:ConvergenceAnalysis:VarianceRelation_Original},
we obtain:
	\begin{align}
		\label{Equ:ConvergenceAnalysis:Adaptation_bound}
		\!\!\mathbb{E}\|\tilde{\bm{\psi}}_{k,i}\|^2	
									&\le	(\gamma_k^2 \!+\! \mu_k^2 \alpha\|S\|_1^2 )
										\!\cdot\!
										\mathbb{E}\|\tilde{\bm{\phi}}_{k,i-1}\|^2
											\!+\!
											\mu_k^2
											\;\|S\|_1^2 \;
											\sigma_v^2
	\end{align}
for $k=1,\ldots,N$.
Finally, introduce the following network mean-square-error vectors (compare with
\eqref{Equ:ConvergenceAnalysis:global_error_vector}):
    \begin{align}
		\mathcal{X}_{i}		=	\begin{bmatrix}
    								\mathbb{E}\|\tilde{\bm\phi}_{1,i}\|^2\\
    								\vdots\\
    								\mathbb{E}\|\tilde{\bm\phi}_{N,i}\|^2
							    \end{bmatrix}, \;
		\mathcal{Y}_{i}		=	\begin{bmatrix}
    								\mathbb{E}\|\tilde{\bm\psi}_{1,i}\|^2\\
    								\vdots\\
    								\mathbb{E}\|\tilde{\bm\psi}_{N,i}\|^2
							    \end{bmatrix}, \;
		\mathcal{W}_{i}		=	\begin{bmatrix}
    								\mathbb{E}\|\tilde{\bm w}_{1,i}\|^2\\
    								\vdots\\
    								\mathbb{E}\|\tilde{\bm w}_{N,i}\|^2
							    \end{bmatrix}
							\nonumber
	\end{align}
and the matrix
	\begin{align}
		\label{Equ:ConvergenceAnalysis:Gamma}
		\Gamma			\;=&\;	\mathrm{diag}
							\left\{
								\gamma_1^2 + \mu_1^2\alpha\|S\|_1^2,
								\;\ldots\;,
								\gamma_N^2 + \mu_N^2\alpha\|S\|_1^2
							\right\}
	\end{align}
Then,
\eqref{Equ:ConvergenceAnalysis:Combination1_bound}%
--\eqref{Equ:ConvergenceAnalysis:Combination2_bound}
 and \eqref{Equ:ConvergenceAnalysis:Adaptation_bound} can be written as
	\begin{align}
		\label{Equ:ConvergenceAnalysis:VarianceRecursion_Inequality_Separate}
		\begin{cases}
			\mathcal{X}_{i-1}	\preceq	P_1^T \mathcal{W}_{i-1}\\
			\mathcal{Y}_i		\preceq	\Gamma \mathcal{X}_{i-1}
									+ \sigma_v^2\|S\|_1^2\Omega^2 \mathds{1}\\
			\mathcal{W}_i		\preceq	P_2^T \mathcal{Y}_i
		\end{cases}
	\end{align}
where the notation $x \preceq y$ denotes that the components of vector $x$ are less than or equal to the corresponding
components of vector $y$. We now recall the following useful fact that for any matrix $F$ with nonnegative entries,
	\begin{align}
		\label{Equ:ConvergenceAnalysis:x_y_Fx_Fy}
		x \preceq y
		\Rightarrow
		Fx \preceq Fy
	\end{align}
This is
because each entry of the vector $Fy-Fx = F(y-x)$ is nonnegative.
Then, combining all three inequalities
in \eqref{Equ:ConvergenceAnalysis:VarianceRecursion_Inequality_Separate} leads to:
	\begin{align}
		\label{Equ:ConvergenceAnalysis:VarianceRecursion_Inequality_Final}
		\boxed{
			\mathcal{W}_i		\preceq	P_2^T \Gamma P_1^T \mathcal{W}_{i-1}
									+
									\sigma_v^2  \|S\|_1^2 \cdot P_2^T  \Omega^2 \mathds{1}
		}
	\end{align}

\subsection{Mean-Square Stability}
Based on \eqref{Equ:ConvergenceAnalysis:VarianceRecursion_Inequality_Final},
we can now prove that, under certain conditions on the step-size parameters $\{\mu_k\}$,
the mean-square-error vector $\mathcal{W}_i$ is bounded as $i\rightarrow \infty$, and
we use this result in the next subsection to evaluate the steady-state MSE for sufficiently small step-sizes.
	\begin{theorem}[Mean-Square Stability]
		\label{Thm:ConvergenceAnalysis:MS_Theorem}
		If the step-sizes $\{\mu_k\}$ satisfy the
		following condition:
			\begin{equation}
				\label{Equ:ConvergenceAnalysis:MS_Theorem:Stepsizes}
				\boxed{
						0	<	\mu_k	<	\min\left\{
													\frac{2\sigma_{k,\max}}
													{\sigma_{k,\max}^2\!+\!\alpha\|S\|_1^2},
													\frac{2\sigma_{k,\min}}
													{\sigma_{k,\min}^2\!+\!\alpha\|S\|_1^2}
												\right\}
				}
			\end{equation}
		for $k=1,\ldots,N$, where $\sigma_{k,\max}$ and $\sigma_{k,\min}$ are defined as
			\begin{align}
				\label{Equ:ConvergenceAnalysis:MS_Theorem:sigma_min_max}
				&\sigma_{k,\max}	\triangleq	\sum_{l=1}^N s_{l,k} \lambda_{l,\max},	\quad
				\sigma_{k,\min}	\triangleq	\sum_{l=1}^N s_{l,k} \lambda_{l,\min}
			\end{align}
		then, as $i\rightarrow \infty$,
			\begin{align}
				\label{Equ:ConvergenceAnalysis:MS_Theorem:W_infty_Bound}
				\boxed{
					\limsup_{i \rightarrow \infty}\| \mathcal{W}_{i}\|_{\infty}	
                                            \le	\frac{
														\displaystyle
														\big(
															\max_{1\le k \le N} \mu_k^2
														\big)
														\cdot
														\|S\|_1^2 \sigma_v^2
													}
													{
														\displaystyle
														1
														-
														\max_{1 \le k \le N}
														(\gamma_k^2 + \mu_k^2 \alpha\|S\|_1^2)
													}
				}
			\end{align}
		where $\| x \|_{\infty}$ denotes the maximum absolute entry of vector $x$.
	\end{theorem}
	\begin{IEEEproof}
		See Appendix \ref{Appendix:Proof_MeanSquaredStability}.
	\end{IEEEproof}
	\vspace{0.5em}

{\color{black}	
If we let $\alpha\!\!=\!\!0$ and $\sigma_v^2\!\!=\!\!0$ in Theorem \ref{Thm:ConvergenceAnalysis:MS_Theorem},
and examine the arguments leading to it, we conclude the validity of the following result,
 which establishes the convergence of the diffusion strategies
\eqref{Equ:DiffusionAdaptation:ATC}--\eqref{Equ:DiffusionAdaptation:CTA}
in the \emph{absence of gradient noise} (i.e., using the true gradient rather than
stochastic gradient --- see \eqref{Equ:DiffusionAdaptation:ATC0}
and \eqref{Equ:DiffusionAdaptation:CTA0}).
}
	\begin{theorem}[Convergence in Noise-free Case]
		\label{Corollary:ConvergenceAnalysis:Convergence_NoiseFreeCase}
		If there is no gradient noise, i.e., $\alpha=0$ and $\sigma_v^2=0$, then
		the mean-square-error vector becomes the deterministic vector
		$\mathcal{W}_i	=	\mathrm{col}\{\|\tilde{ w}_{1,i}\|^2,\cdots,\|\tilde{ w}_{N,i}\|^2\}$,
		and its entries converge to zero if the step-sizes $\{\mu_k\}$ satisfy the following condition:
			\begin{equation}
				\label{Equ:ConvergenceAnalysis:NoiseFree_Corollary:Stepsizes}
				\boxed{
						0	<	\mu_k	<		\frac{2}{\sigma_{k,\max}}
				}
			\end{equation}
		for $k=1,\ldots,N$, where $\sigma_{k,\max}$
		was defined in \eqref{Equ:ConvergenceAnalysis:MS_Theorem:sigma_min_max}.	
		\hfill\IEEEQED
	\end{theorem}
	
We observe that, in the absence of noise, the deterministic error vectors,
$\tilde{w}_{k,i}$, will tend to zero as $i\rightarrow \infty$ even with constant
(i.e., non-vanishing) step-sizes. This result implies the interesting fact that, in the noise-free case,
the nodes can reach agreement
\emph{without} the need to impose diminishing step-sizes.

\subsection{Steady-State Performance}
\label{Sec:ConvergenceAnalysis:SteadyState}
Expression \eqref{Equ:ConvergenceAnalysis:MS_Theorem:Stepsizes}
provides a condition on the step-size parameters $\{\mu_k\}$ to ensure the
mean-square stability of the diffusion strategies
\eqref{Equ:DiffusionAdaptation:ATC}--\eqref{Equ:DiffusionAdaptation:CTA}. At the same time, expression
\eqref{Equ:ConvergenceAnalysis:MS_Theorem:W_infty_Bound} gives an upper bound on
how large $\mathcal{W}_i$ can be at steady-state.
Since the $\infty$-norm of a vector is defined as the largest absolute value of its entries,  then
\eqref{Equ:ConvergenceAnalysis:MS_Theorem:W_infty_Bound} bounds
the MSE of the worst-performing node in the network.
We can derive closed-form expressions for MSEs when the step-sizes are assumed
to be sufficiently small.
Indeed, we first conclude from \eqref{Equ:ConvergenceAnalysis:MS_Theorem:W_infty_Bound}
that for step-sizes that are sufficiently small,
each $\bm{w}_{k,i}$ will get closer to $w^o$ at steady-state.
To verify this fact, assume the step-sizes are small enough so that the nonnegative factor
$\gamma_k$ that was defined earlier in
\eqref{Equ:ConvergenceAnalysis:gamma_k} becomes
	\begin{align}
		\label{Equ:ConvergenceAnalysis:gamma_k_SmallStepsize}
		\gamma_k	=	1 - \mu_k \sum_{l=1}^N s_{l,k} \lambda_{l,\min}
				=	1 - \mu_k \sigma_{k,\min}
	\end{align}
where $\sigma_{k,\min}$ was given by \eqref{Equ:ConvergenceAnalysis:MS_Theorem:sigma_min_max}.
Substituting \eqref{Equ:ConvergenceAnalysis:gamma_k_SmallStepsize} into
\eqref{Equ:ConvergenceAnalysis:MS_Theorem:W_infty_Bound}, we obtain:
	\begin{align}
		&\limsup_{i\rightarrow\infty}
        \|\mathcal{W}_{i}\|_{\infty}	
                                        \nonumber\\
                                &\quad
                                    \le	\frac{
											\displaystyle
											\Big(\max_{1\le k \le N} \mu_k^2\Big)
											\cdot
											\|S\|_1^2 \sigma_v^2
										}
										{
											\displaystyle
											1
											\!-\!
											\max_{1 \le k \le N}
                                            \Big\{
											        (1\!-\!\mu_k\sigma_{k,\min})^2
                                                    \!+\! \mu_k^2 \alpha \|S\|_1^2
                                            \Big\}
										}
                                        \nonumber\\
								&\quad
                                    \le	\frac{  \displaystyle\Big(\max_{1\le k \le N} \mu_k^2\Big)
                                                \cdot
                                                \|S\|_1^2 \sigma_v^2
                                        }
                                        {\displaystyle
                                            \min_{1\le k \le N}
                                            \Big\{
                                                    \mu_{k}
                                                    \Big[
                                                        2\sigma_{k,\min}
                                                        -
                                                        \mu_k
                                                        (\sigma_{k,\min}^2+\alpha \|S\|_1^2)
                                                    \Big]
                                            \Big\}
                                        }
                                        \nonumber\\
		\label{Equ:ConvergenceAnalysis:MSD_Vanish_SmallStepSize}
                                &\quad
                                    \le
                                    \frac{
                                            \|S\|_1^2 \sigma_v^2
                                        }
                                        {\displaystyle
                                            \min_{1\le k \le N}
                                            \Big\{
                                                        2\sigma_{k,\min}
                                                        \!-\!
                                                        \mu_k
                                                        (\sigma_{k,\min}^2\!+\!\alpha \|S\|_1^2)
                                            \Big\}
                                        }
                                    \cdot
                                    \frac{\mu_{\max}^2}{\mu_{\min}}
	\end{align}
where
	\begin{align}
        \label{Equ:PerformanceAnalysis:mu_max_min_def}
		\mu_{\max}	\!\triangleq 	\!\displaystyle\max_{1\le k\le N} \mu_k,
		\qquad
		\mu_{\min}	\!\triangleq 	\!\displaystyle\min_{1\le k\le N}\mu_k
	\end{align}
For sufficiently small step-sizes, the denominator in \eqref{Equ:ConvergenceAnalysis:MSD_Vanish_SmallStepSize}
can be approximated as
    \begin{align}
        \label{Equ:PerformanceAnalysis:W_inf_bound_intermediate_approx}
        2\sigma_{k,\min} \!-\! \mu_k (\sigma_{k,\min}^2\!+\!\alpha \|S\|_1^2)
            \approx 2\sigma_{k,\min}
    \end{align}
Substituting
into \eqref{Equ:ConvergenceAnalysis:MSD_Vanish_SmallStepSize}, we get
	\begin{align}
		\label{Equ:ConvergenceAnalysis:MSD_Vanish_SmallStepSize_final}
        \boxed{
		\limsup_{i\rightarrow\infty}
        \|{\mathcal{W}}_{i}\|_{\infty}		\;\le\;		\frac{\|S\|_1^2 \sigma_v^2}
												{\displaystyle2\min_{1\le k \le N} \sigma_{k,\min}}
											\cdot \frac{\mu_{\max}^2}{\mu_{\min}}
        }
	\end{align}
Therefore, if the step-sizes are sufficiently small,
the MSE of each node
becomes small as well.
This result is clear when all nodes
use the same step-sizes such that $\mu_{\max}=\mu_{\min}=\mu$.
Then, the right-hand side of \eqref{Equ:ConvergenceAnalysis:MSD_Vanish_SmallStepSize_final}
is on the order of $O(\mu)$, as indicated.
It follows that $\{\tilde{\bm{w}}_{k,i}\}$ are small in the mean-square-error sense
at small step-sizes, which also means that
the mean-square value of $\tilde{\bm{\phi}}_{k,i-1}$ is small because it is a convex combination
of $\{\tilde{\bm{w}}_{k,i}\}$ (recall \eqref{Equ:ConvergenceAnalysis:phi_ki}).
Then, by definition \eqref{Equ:ConvergenceAnalysis:H_kim1},
{\color{black}
in steady-state (for large enough $i$),
}
the matrix $\bm{H}_{l,k,i-1}$
can be approximated by:
	\begin{align}
        \label{Equ:PerformanceAnalysis:H_lkim1_approx_steadystate}
		\bm{H}_{l,k,i-1}	\approx	\int_{0}^{1}
							\nabla^2
							J_l(w^o) dt
					=		\nabla^2
							J_l(w^o)
	\end{align}
In this case, the matrix $\bm{H}_{l,k,i-1}$ is not random anymore and is
not dependent on the error vector $\tilde{\bm{\phi}}_{k,,i-1}$. Accordingly,
in steady-state,
the matrix
$\bm{\mathcal{D}}_{i-1}$ that was defined in \eqref{Equ:ConvergenceAnalysis:D_i_minus_1}
is not random anymore and it becomes
	\begin{align}
		\label{Equ:ConvergenceAnalysis:D_inf}
		\boxed{
		\bm{\mathcal{D}}_{i-1}	
					\!\approx\!
						\mathcal{D}_{\infty}
					\!\triangleq\!	
						\sum_{l=1}^N
						\mathrm{diag}
						\big\{
							 s_{l,1} \nabla_w^2
							J_l(w^o),
							\cdots,\!
							s_{l,N} \nabla_w^2
							J_l(w^o)
						\big\}
		}
	\end{align}
{\color{black}
As a result,
in steady-state,
the original error
recursion \eqref{Equ:ConvergenceAnalysis:ErrorRecursion_final} can be approximated
by
    \begin{align}
		\label{Equ:ConvergenceAnalysis:ErrorRecursion_final_approx}
		\boxed{
			\tilde{\bm w}_i	=	\mathcal{P}_2^T
							[I_{MN}-\mathcal{M}{\mathcal{D}}_{\infty}]
							\mathcal{P}_1^T
							\tilde{\bm w}_{i-1}
							+
							\mathcal{P}_2^T \mathcal{M} \bm{g}_i
		}
	\end{align}
Taking expectations of both sides of \eqref{Equ:ConvergenceAnalysis:ErrorRecursion_final_approx},
we obtain the following mean-error recursion
    \begin{align}
        \label{Equ:PerformanceAnalysis:MeanErrorRecursion_Approx}
        \E\tilde{\bm w}_i	
                        =	\mathcal{P}_2^T
							[I_{MN}-\mathcal{M}{\mathcal{D}}_{\infty}]
							\mathcal{P}_1^T
                            \cdot
							\E\tilde{\bm w}_{i-1}, \quad i \rightarrow\infty
    \end{align}
which converges to zero if the matrix
    \begin{align}
        \label{Equ:PerformanceAnalysis:B_cal_def}
        \mathcal{B} \triangleq \mathcal{P}_2^T[I_{MN}-\mathcal{M}{\mathcal{D}}_{\infty}]\mathcal{P}_1^T
    \end{align}
is stable.
The stability of $\mathcal{B}$ can be guaranteed when the step-sizes are
sufficiently small (or chosen according to \eqref{Equ:ConvergenceAnalysis:MS_Theorem:Stepsizes})
--- see the proof in Appendix \ref{Appendix:Stability_F}. Therefore, in steady-state, we have
    \begin{align}
        \label{Equ:PerformanceAnalysis:Ew_inf_zero}
        \boxed{
            \lim_{i \rightarrow \infty} \E\tilde{\bm{w}}_i = 0
        }
    \end{align}
}%
Next, we determine an expression (rather than a bound) for the MSE.
To do this, we need to evaluate the covariance matrix of the gradient noise vector $\bm{g}_i$.
Recall from \eqref{Equ:ConvergenceAnalysis:G_i} that $\bm{g}_i$ depends
on $\{{\bm{\phi}}_{k,i-1}\}$, which is close to $w^o$ at steady-state for small step-sizes.
Therefore, it is sufficient to determine the covariance matrix of $\bm{g}_i$ at $w^o$.
We denote this covariance matrix by:
			\begin{align}
				R_v         \;\triangleq\;& \mathbb{E}\{\bm{g}_i\bm{g}_i^T\}\big|_{\phi_{k,i-1}=w^o}
                                            \nonumber\\
                            \;=\;&          \mathbb{E}
                                            \bigg\{\!
                                            \Big[\!
                                                \sum_{l=1}^N
						                          \mathrm{col}
						                          \big\{
							                             s_{l,1}\bm{v}_{l,i}(w^o),
							                             \cdots,
							                             s_{l,N}\bm{v}_{l,i}(w^o)
						                          \big\}\!
                                            \Big]
                                            \nonumber\\
				\label{Equ:ConvergenceAnalysis:GradientNoise_PreciseModel}
                                            &\times
                                            \Big[\!
                                                \sum_{l=1}^N
						                          \mathrm{col}
						                          \big\{
							                             s_{l,1}\bm{v}_{l,i}(w^o),
							                             \cdots,
							                             s_{l,N}\bm{v}_{l,i}(w^o)
						                          \big\}\!
                                            \Big]^T\!
                                            \bigg\}
			\end{align}
In practice, we can evaluate $R_v$ from the expressions
of $\{\bm{v}_{l,i}(w^o)\}$. For example, for the case of the quadratic
cost \eqref{Equ:PerformanceAnalysis:QuadraticCost_LMS}, we can
substitute \eqref{Equ:PerformanceAnalysis:GradientNoise_LMS}
into \eqref{Equ:ConvergenceAnalysis:GradientNoise_PreciseModel}
to evaluate $R_v$.

Returning to the last term in the first equation of
\eqref{Equ:ConvergenceAnalysis:WeightedEnergyConservation_Relation},
we can evaluate it as follows:
	\begin{align}
		\mathbb{E}\|\mathcal{P}_2^T \mathcal{M} \bm{g}_i\|_\Sigma^2
				\;=&\;
					\mathbb{E}
					\bm{g}_i^T\mathcal{M} \mathcal{P}_2
					\Sigma
					\mathcal{P}_2^T \mathcal{M} \bm{g}_i
                    \nonumber\\
				\;=&\;	
					\mathrm{Tr}
					\left(
						\Sigma	
						\mathcal{P}_2^T \mathcal{M}
						\mathbb{E}\{
						\bm{g}_i
						\bm{g}_i^T
						\}
						\mathcal{M} \mathcal{P}_2				
					\right)
					\nonumber\\
				\;=&\;
		\label{Equ:ConvergenceAnalysis:NoiseTerm_PreciseModel}
					\mathrm{Tr}
					\left(
						\Sigma	
						\mathcal{P}_2^T \mathcal{M}
						R_v
						\mathcal{M} \mathcal{P}_2				
					\right)
	\end{align}
Using \eqref{Equ:ConvergenceAnalysis:D_inf}, the matrix $\bm{\Sigma}'$ in
\eqref{Equ:ConvergenceAnalysis:WeightedEnergyConservation_Relation} becomes a
deterministic quantity as well, and is given by:
	\begin{align}
		\label{Equ:ConvergenceAnalysis:Sigma_Prime_Approx_SmallStepSize}
		{\Sigma}'	\approx	\mathcal{P}_1
							[I_{MN}-\mathcal{M}{\mathcal{D}}_{\infty}]
							\mathcal{P}_2
							\Sigma
							\mathcal{P}_2^T
							[I_{MN}-\mathcal{M}{\mathcal{D}}_{\infty}]
							\mathcal{P}_1^T
	\end{align}
Substituting \eqref{Equ:ConvergenceAnalysis:NoiseTerm_PreciseModel}
and \eqref{Equ:ConvergenceAnalysis:Sigma_Prime_Approx_SmallStepSize} into
\eqref{Equ:ConvergenceAnalysis:WeightedEnergyConservation_Relation},
an approximate variance relation is obtained for small step-sizes:
	\begin{align}
		\label{Equ:ConvergenceAnalysis:WeightedEnergyConservation_Relation_Approx_SmallStepSize}
			\!\!\!\!\mathbb{E}\|\tilde{\bm{w}}_{i}\|_\Sigma^2
					&\approx
							\mathbb{E}\|\tilde{\bm{w}}_{i-1}\|_{{\Sigma}'}^2
							+\mathrm{Tr}
							\left(
								\Sigma
								\mathcal{P}_2^T \mathcal{M}
								R_v
								\mathcal{M} \mathcal{P}_2					
							\right)
							\\
		\label{Equ:ConvergenceAnalysis:Sigma_PrimePrime}
			\!\!\!\!{\Sigma}'	&\approx
							\mathcal{P}_1
							[I_{MN}\!\!-\!\!\mathcal{M}{\mathcal{D}}_{\infty}]
							\mathcal{P}_2
							\Sigma
							\mathcal{P}_2^T
							[I_{MN}\!\!-\!\!\mathcal{M}{\mathcal{D}}_{\infty}]
							\mathcal{P}_1^T
	\end{align}	
Let $\sigma=\mathrm{vec}(\Sigma)$ denote the vectorization operation that
stacks the columns of a matrix $\Sigma$ on top of each other. We shall use
the notation $\|x\|_{\sigma}^2$ and $\|x\|_{\Sigma}^2$
interchangeably to denote the weighted squared Euclidean norm of a vector.
Using the Kronecker product property\cite[p.147]{laub2005matrix}:
		$\mathrm{vec}(U\Sigma V) = (V^T \otimes U ) \mathrm{vec}(\Sigma)$,
we can vectorize ${\Sigma}'$ in \eqref{Equ:ConvergenceAnalysis:Sigma_PrimePrime}
and find that its vector form is related to $\Sigma$ via the following \emph{linear} relation:
		$\sigma'	\triangleq	\mathrm{vec}(\Sigma')	\approx	\mc{F}\sigma$,
where, for sufficiently small steps-sizes (so that higher powers of the step-sizes can be ignored), the
matrix $\mc{F}$ is given by
	\begin{align}
		\label{Equ:ConvergenceAnalysis:F}
		\!\!\!\boxed{
		\mc{F}	
            \!\triangleq\!
				\big(
					\mathcal{P}_1
					[I_{MN}\!-\!\mathcal{M}{\mathcal{D}}_{\infty}]
					\mathcal{P}_2
				\big)\!
				\otimes\!
				\big(
					\mathcal{P}_1
					[I_{MN}\!-\!\mathcal{M}{\mathcal{D}}_{\infty}]
					\mathcal{P}_2
				\big)
		}
	\end{align}
Here, we used the fact that $\mathcal{M}$ and $\mathcal{D}_\infty$ are block diagonal
and symmetric.
Furthermore, using the property $\mathrm{Tr}(\Sigma X) = \mathrm{vec}(X^T)^T \sigma$,
we can rewrite \eqref{Equ:ConvergenceAnalysis:WeightedEnergyConservation_Relation_Approx_SmallStepSize}
as
	\begin{align}
		\label{Equ:ConvergenceAnalysis:WeightedEnergyConservation_final}
		\mathbb{E}\|\tilde{\bm{w}}_{i}\|_\sigma^2
					\;\approx\;&	
							\mathbb{E}\|\tilde{\bm{w}}_{i-1}\|_{\mc{F}\sigma}^2
							+
							\left[
								\mathrm{vec}
								\left(
									\mathcal{P}_2^T \mathcal{M}
									R_v
									\mathcal{M} \mathcal{P}_2	
								\right)
							\right]^T \!\!
							\sigma
	\end{align}
It is shown in \cite[pp.344--346]{Sayed08} that recursion \eqref{Equ:ConvergenceAnalysis:WeightedEnergyConservation_final} converges
to a \emph{steady-state value} if the matrix $\mc{F}$ is stable.
This condition is guaranteed
when the step-sizes are sufficiently small (or chosen according to
\eqref{Equ:ConvergenceAnalysis:MS_Theorem:Stepsizes}) --- see Appendix \ref{Appendix:Stability_F}.
Finally, denoting
    \begin{align}
        \mathbb{E}\|\tilde{\bm{w}}_{\infty}\|_{\sigma}^2 \triangleq \lim_{i\rightarrow\infty}
                    \mathbb{E}\|\tilde{\bm{w}}_{i}\|_{\sigma}^2
    \end{align}
and letting $i \rightarrow \infty$, expression \eqref{Equ:ConvergenceAnalysis:WeightedEnergyConservation_final}
becomes
	\begin{align}
		\mathbb{E}\|\tilde{\bm{w}}_{\infty}\|_\sigma^2
					\;\approx\;&	
							\mathbb{E}\|\tilde{\bm{w}}_{\infty}\|_{\mc{F}\sigma}^2
							+
							\left[
								\mathrm{vec}
								\left(
									\mathcal{P}_2^T \mathcal{M}
									R_v
									\mathcal{M} \mathcal{P}_2	
								\right)
							\right]^T
							\sigma
							\nonumber
	\end{align}
so that
	\begin{align}
		\label{Equ:ConvergenceAnalysis:SteadyStatePerformance_final}
		\boxed{
		\mathbb{E}\|\tilde{\bm{w}}_{\infty}\|_{(I-\mc{F})\sigma}^2
				\approx	\left[
						\mathrm{vec}
						\left(
							\mathcal{P}_2^T \mathcal{M}
							R_v
							\mathcal{M} \mathcal{P}_2	
						\right)
					\right]^T
					\sigma
		}
	\end{align}
Expression \eqref{Equ:ConvergenceAnalysis:SteadyStatePerformance_final}
is a useful result: it allows us to derive several performance
metrics through the proper selection of the free weighting parameter $\sigma$ (or $\Sigma$).
First, to be able to evaluate steady-state performance metrics from
\eqref{Equ:ConvergenceAnalysis:SteadyStatePerformance_final},
we need $(I-\mc{F})$ to be invertible, which is guaranteed by the stability of matrix $\mc{F}$ --- see Appendix \ref{Appendix:Stability_F}.
Given that $(I-\mc{F})$ is a stable matrix, we can now resort to
\eqref{Equ:ConvergenceAnalysis:SteadyStatePerformance_final} and use it to
evaluate various performance metrics by choosing proper weighting matrices
$\Sigma$ (or $\sigma$),
as it was done in \cite{Cattivelli10} for the mean-square-error estimation problem. For example,
the MSE of any node $k$ can be obtained by computing $\mathbb{E}\|\tilde{\bm{w}}_{\infty}\|_{T}^2$
with a block weighting matrix $T$ that has an identity matrix at block $(k,k)$ and zeros elsewhere:
	\begin{align}
		\mathbb{E}\|\tilde{\bm{w}}_{k,\infty}\|^2	=	\mathbb{E}\|\tilde{\bm{w}}_{\infty}\|_{T}^2
	\end{align}
Denote the vectorized version of this matrix by
$t_k$, i.e.,
	\begin{align}
		t_k	\triangleq	\mathrm{vec}(\mathrm{diag}(e_k)\otimes I_M)
	\end{align}
where $e_k$ is a vector whose $k$th entry is one and zeros elsewhere.
Then, if we select $\sigma$ in \eqref{Equ:ConvergenceAnalysis:SteadyStatePerformance_final}
as $\sigma = (I-\mc{F})^{-1}t_k$,
the term on the left-hand side becomes the desired $\mathbb{E}\|\tilde{\bm{w}}_{k,\infty}\|^2$
and MSE for node $k$  is therefore given by:
	\begin{align}
		\label{Equ:ConvergenceAnalysis:MSD_k}
			\mathrm{MSE}_k	\approx	
								\left[
									\mathrm{vec}
									\left(
										\mathcal{P}_2^T \mathcal{M}
										R_v
										\mathcal{M} \mathcal{P}_2	
									\right)
								\right]^T
								(I-\mc{F})^{-1} t_k
	\end{align}
This value for $\mathrm{MSE}_k$ is actually the $k$th entry of $\mathcal{W}_{\infty}$
defined as
    \begin{align}
        \mathcal{W}_\infty \triangleq \lim_{i\rightarrow\infty} \mathcal{W}_i
    \end{align}
Then, we arrive at an expression for $\mathcal{W}_{\infty}$
(as opposed to the bound for it in \eqref{Equ:ConvergenceAnalysis:MS_Theorem:W_infty_Bound},
as was explained earlier;
expression \eqref{Equ:ConvergenceAnalysis:W_infty} is derived under the assumption of sufficiently small step-sizes):
	\begin{align}
		\label{Equ:ConvergenceAnalysis:W_infty}
		\!\!\boxed{
		\mathcal{W}_{\infty}	\!\approx\!
							\left\{\!
								I_N \!\otimes\!
								\left(
									\left[
										\mathrm{vec}
										\left(
											\mathcal{P}_2^T \mathcal{M}
											R_v
											\mathcal{M} \mathcal{P}_2	
										\right)
									\right]^T
									(I\!\!-\!\!\mc{F})^{-1} 								
								\right)\!
							\right\}
							t
		}
	\end{align}
where $t=\mathrm{col}\{t_1,\ldots,t_N\}$.
If we are interested in the network MSE, then the weighting matrix of
$\mathbb{E}\|\tilde{\bm{w}}_{\infty}\|_{T}^2$ should be chosen as $T=I_{MN}/N$.
Let $q$ denote the vectorized version of $I_{MN}$, i.e.,
	\begin{align}
		q \triangleq \mathrm{vec}(I_{MN})
	\end{align}
and select $\sigma$ in \eqref{Equ:ConvergenceAnalysis:SteadyStatePerformance_final}
as $\sigma = (I\!-\!\mc{F})^{-1} q/N$.
The network MSE is then given by
	\begin{equation}
		\label{Equ:ConvergenceAnalysis:MSD_Network}
		\boxed{
			\begin{split}
			\overline{\mathrm{MSE}}
				&\triangleq	\frac{1}{N} \sum_{k=1}^N \mathrm{MSE}_k\\
                &\approx		\frac{1}{N}
							\left[
								\mathrm{vec}
								\left(
									\mathcal{P}_2^T \mathcal{M}
									R_v
									\mathcal{M} \mathcal{P}_2	
								\right)
							\right]^T
							(I-\mc{F})^{-1} q		
			\end{split}
		}
	\end{equation}
{\color{black}
The approximate expressions \eqref{Equ:ConvergenceAnalysis:W_infty} and
\eqref{Equ:ConvergenceAnalysis:MSD_Network} hold when the step-sizes are
small enough so that \eqref{Equ:ConvergenceAnalysis:D_inf} holds. In the next
section, we will
see that they are consistent with the simulation results.
}

\section{Simulation Results}
\label{Sec:Simulation}
In this section we illustrate the performance of the diffusion strategies
\eqref{Equ:DiffusionAdaptation:ATC}--\eqref{Equ:DiffusionAdaptation:CTA}
by considering two applications.
We consider a randomly generated connected network topology with a cyclic path.
There are a total of $N=10$ nodes in the network, and nodes are assumed connected when they are
close enough geographically.
In the simulations, we consider two applications: a regularized least-mean-squares estimation problem
with sparse parameters,
and a collaborative localization problem.

\subsection{Distributed Estimation with Sparse Data}
\label{Sec:Simulation:L1}

\begin{figure*}[t!]
   \centerline{
		   \subfigure[Learning curves ($\gamma=2$ and $\epsilon=10^{-3}$).]
		   {
			   \includegraphics[width=3in]%
			   {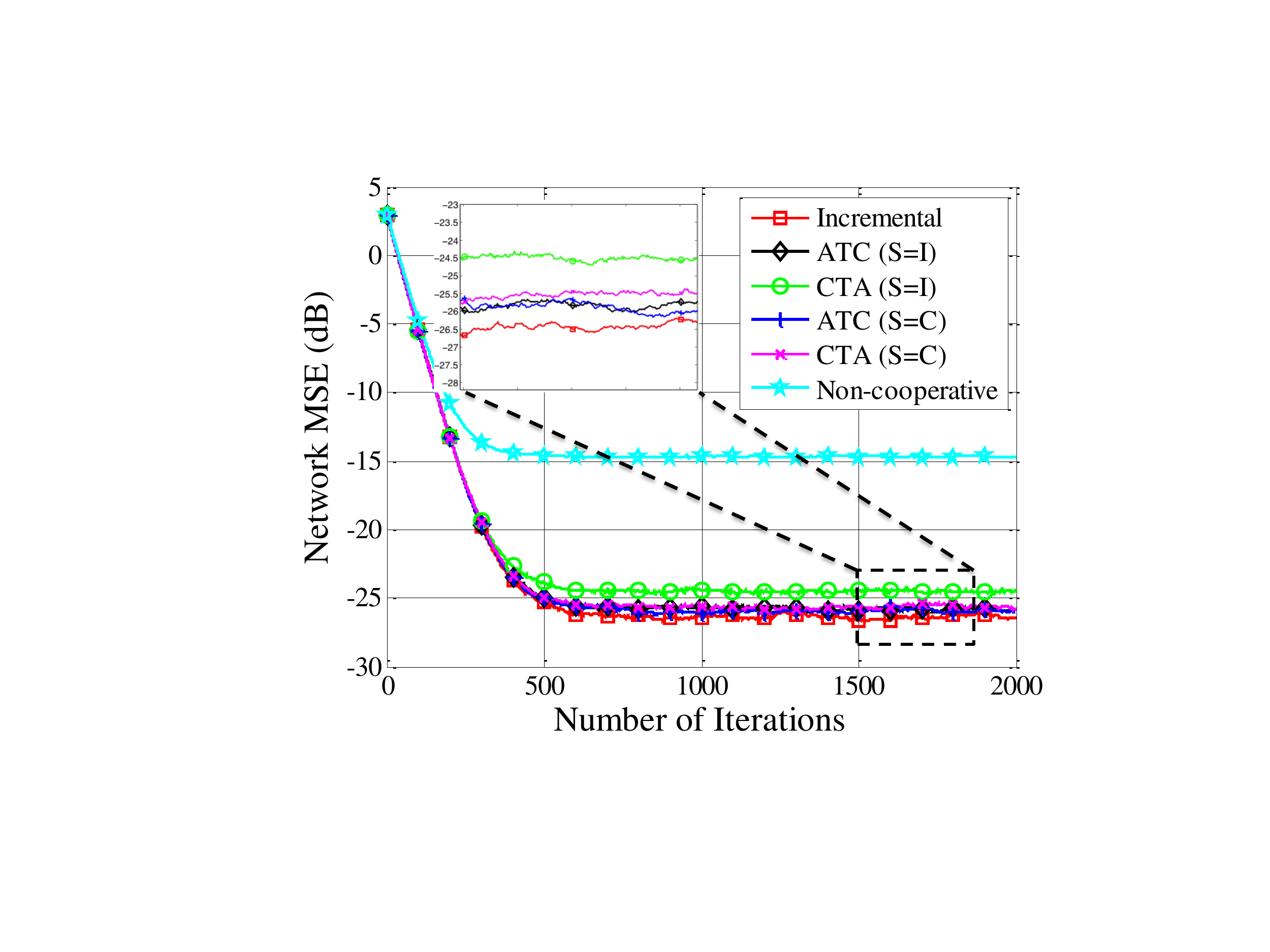}
			   \label{fig:Simulation:LearningCurve}
		   }
    \hfil
   \subfigure[Steady-state MSD ($\mu=10^{-3}$).]
  		{
				   \includegraphics[width=3.1in]%
				   {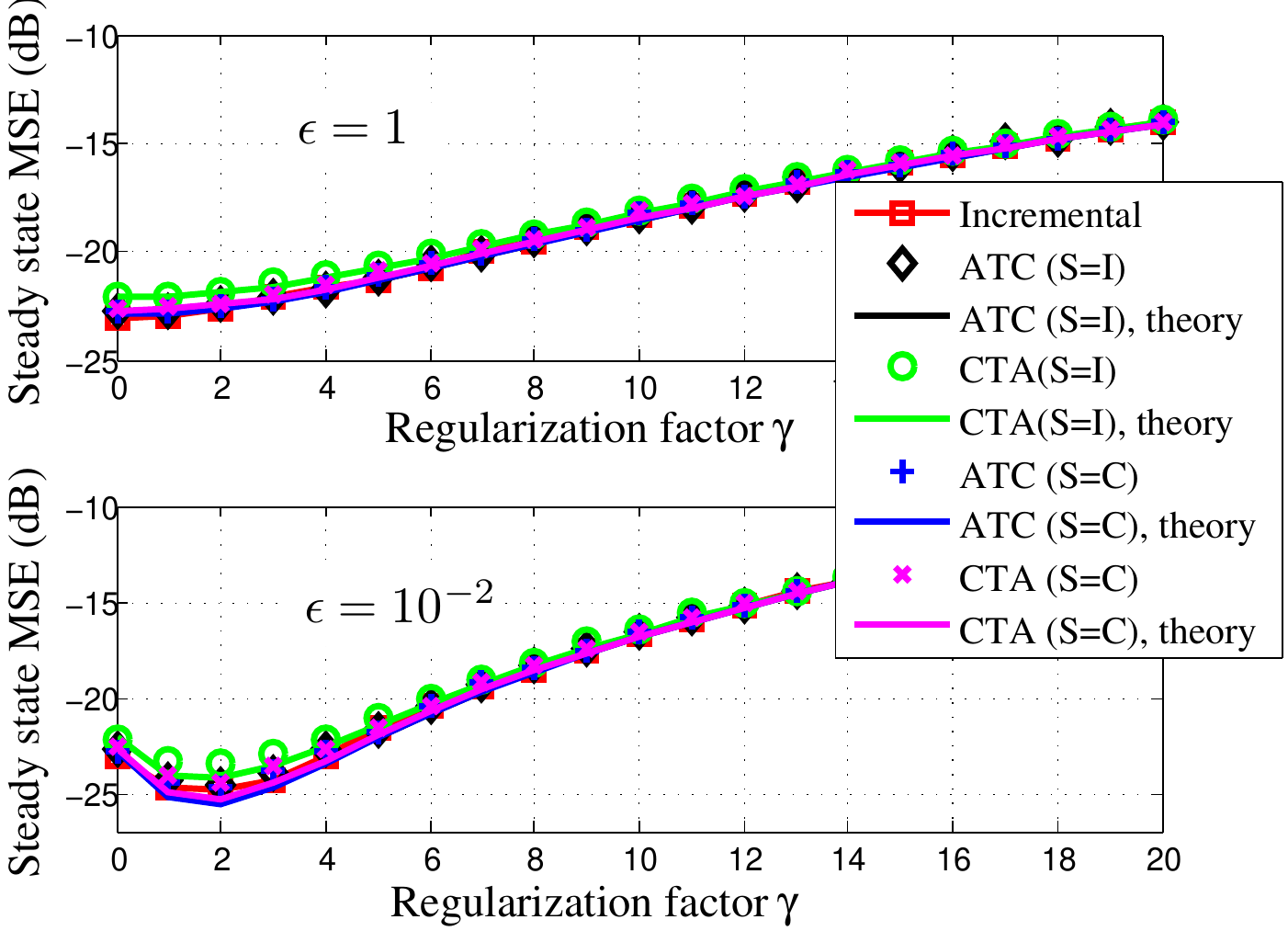}
				   \label{fig:Simulation:MSDvsGamma}
 			    }
			   }
\caption{Transient and steady-state performance of distributed estimation with sparse parameters.}
\end{figure*}

Assume each node $k$ has access to data $\{ \bm{U}_{k,i},  \bm{d}_{k,i}\}$, generated according to
the following model:
	\begin{align}
		\label{Equ:Simulation:RegularizedLMS:LinearModel}
		 \bm{d}_{k,i}	=	 \bm{U}_{k,i} {w}^o +  \bm{z}_{k,i}
	\end{align}
where $\{\bm{U}_{k,i}\}$ is a  sequence of $K \times M$ i.i.d. Gaussian random  matrices.
The entries of each $\bm{U}_{k,i}$ have zero mean
and unit variance,
and $ \bm{z}_{k,i} \sim \mathcal{N}(0,\sigma_z^2I_K)$ is the measurement noise that is
temporally and spatially white and is
independent of $ \bm{U}_{l,j}$ for all $k,l,i,j$.
Our objective is to estimate ${w}^o$ from the data set $\{ \bm{U}_{k,i},  \bm{d}_{k,i}\}$ in a distributed manner.
In many applications, the vector ${w}^o$ is sparse such as
$$
		{w}^o	=	[
					1 \; 0	 \; \ldots \; 0 \; 1
					]^T \in \mathbb{R}^M
$$
One way to search for sparse solutions is to consider a global cost function of the following form
\cite{diLorenzo2012icassp}:
	\begin{align}
		\label{Equ:Simulation:J_glob_L1RLS}
		J^{\mathrm{glob}}(w)		=	\sum_{l=1}^N
								\mathbb{E}\| \bm{d}_{l,i} - \bm{U}_{l,i} w \|_2^2 + \rho R(w)
	\end{align}
where $R(w)$ and $\rho$ are the regularization function and regularization factor,  respectively.
A popular choice is $R(w)=\|w\|_1$, which helps enforce sparsity and is convex
\cite{tibshirani1996regression,baraniuk2007compressive,candes2008enhancing,mateos2010distributed,
kopsinis2011online,
diLorenzo2012icassp}.
However, this choice is non-differentiable, and we would need to
apply sub-gradient methods \cite[pp.138--144]{poliak1987introduction} for a proper implementation.
Instead, we use the following twice-differentiable approximation for $\|w\|_1$:
	\begin{align}
		\label{Equ:Simulation:R_w}
		R(w)	=	\sum_{m=1}^M \sqrt{[w]_m^2 + \epsilon^2}
	\end{align}
where $[w]_m$ denotes the $m$-th entry of $w$, and $\epsilon$ is a small number.
We see that, as $\epsilon$ goes to zero, $R(w) \approx \|w\|_1$.
Obviously, $R(w)$ is convex,
and we can apply the diffusion algorithms
to minimize  \eqref{Equ:Simulation:J_glob_L1RLS} in a distributed manner.
To do so, we decompose the global cost into a sum of $N$ individual costs:
	\begin{align}
		\label{Equ:Simulation:J_k_L1RLS}
		J_l(w)	=	\mathbb{E}\| \bm{d}_{l,i}  -  \bm{U}_{l,i} w \|_2^2 + \frac{\rho}{N} R(w)		
	\end{align}
for $l=1,\ldots,N$.
Then, using algorithms \eqref{Equ:DiffusionAdaptation:ATC0}
and \eqref{Equ:DiffusionAdaptation:CTA0}, each node $k$ would update its estimate of ${w}^o$ by using
the gradient vectors of $\{J_l(w)\}_{l \in \mathcal{N}_k}$, which are given by:
	\begin{align}
		\nabla_w J_l(w)		\;=\;&	2\mathbb{E}\left(\bm{U}_{l,i}^T\bm{U}_{l,i}\right) {w}
							-
							2\mathbb{E}\left( \bm{U}_{l,i}^T \bm{d}_{l,i}\right)
                            \nonumber\\
							\;&+
							\frac{\rho}{N} \nabla_w R(w)
	\end{align}
However, the nodes are assumed to have access to measurements $\{{U}_{l,i},{d}_{l,k}\}$
and not to the second-order moments $\mathbb{E}\big(\bm{U}_{l,i}^T \bm{U}_{l,i}\big)$ and
$\mathbb{E}\big(\bm{U}_{l,i}^T\bm{d}_{l,i}\big)$. In this case, nodes can use the available
measurements to approximate the gradient vectors in
\eqref{Equ:DiffusionAdaptation:ATC}
and \eqref{Equ:DiffusionAdaptation:CTA} as:
	\begin{align}
			\widehat{\nabla}_w J_l(w)
				=	2{U}_{l,i}^T
					\left[	
						{U}_{l,i} w \!-\! {d}_{l,i}
					\right]
					\!+\!
					\frac{\rho}{N}
					\nabla_w R(w)
	\end{align}
where
	\begin{align}
		\nabla_w R(w)	=	\begin{bmatrix}
							\displaystyle\frac{[w]_1}{\sqrt{[w]_1^2 + \epsilon^2}}
							&
							\cdots
							&
							\displaystyle\frac{[w]_M}{\sqrt{[w]_M^2 + \epsilon^2}}
						\end{bmatrix}^T
	\end{align}
In the simulation, we set $M=50$, $K=5$, $\sigma_v^2=1$,
and ${w}^o=[1 \; 0	 \; \ldots \; 0 \; 1]^T$.
We apply both diffusion and incremental methods to solve the distributed learning problem,
where the incremental approach
\cite{bertsekas1997new, nedic2001incremental, rabbat2005quantized, lopes2007incremental}
uses the following construction to determine $\bm{w}_i$:
    	\begin{itemize}
            \item
    			Start with $\bm{\psi}_{0,i}		=	\bm{w}_{i-1}$
    			at the node at the beginning of the incremental cycle.
            \item
			    Cycle through the nodes $k=1,\ldots,N$:
                \begin{align}
			         \bm{\psi}_{k,i}	=	\bm{\psi}_{k-1,i}
                                            - \mu\widehat{\nabla}_w J_{k}(\bm{\psi}_{k-1,i})
                \end{align}
            \item
			 Set $\bm{w}_i			\leftarrow	\bm{\psi}_{N,i}$.
            \item
			 Repeat.
		\end{itemize}
The results are averaged over $100$ trials.
The step-sizes for ATC, CTA and non-cooperative algorithms are set to $\mu=10^{-3}$,
and the step-size for the incremental algorithm is set to $\mu=10^{-3}/N$.
This is because the incremental algorithm cycles through all $N$ nodes every iteration.
We therefore need to ensure the same convergence rate for both algorithms for a fair comparison
\cite{takahashi2010diffusion}. For ATC and CTA strategies, we use simple
averaging weights for the combination step, and for ATC and CTA with
gradient exchange, we use Metropolis weights for $\{c_{l,k}\}$ to combine the gradients
(see Table III in \cite{Cattivelli10} for the definitions of averaging weights and Metropolis weights).
We use expression \eqref{Equ:ConvergenceAnalysis:MSD_Network}
to evaluate the theoretical performance of the diffusion strategies.
{\color{black}
As a remark, expression \eqref{Equ:ConvergenceAnalysis:MSD_Network} gives
the MSE with respect to the minimizer of the cost $J^{\mathrm{glob}}(w)$ in
\eqref{Equ:Simulation:J_glob_L1RLS}. In this example, the minimizer of the cost
\eqref{Equ:Simulation:J_glob_L1RLS}, denoted as $\hat{w}^o$, is biased away from the
model parameter $w^o$
in \eqref{Equ:Simulation:RegularizedLMS:LinearModel} when
the regularization factor $\gamma \neq 0$.
To evaluate the theoretical MSE with respect to $w^o$, we use
    \begin{align}
        \overline{\mathrm{MSD}}
                                        &=   \lim_{i \rightarrow \infty}
                                            \frac{1}{N}
                                            \sum_{k=1}^N
                                            \E
                                            \|
                                                w^o-\bm{w}_{k,i}
                                            \|^2
                                            \nonumber\\
        \label{Equ:Simulation:L1_MSD_modified}
                                        &=  \E
                                            \|
                                                w^o - \hat{w}^o
                                            \|^2
                                            +
                                            \lim_{i \rightarrow \infty}
                                            \frac{1}{N}
                                            \sum_{k=1}^N
                                                    \E
                                                    \|
                                                        \hat{w}^o-\bm{w}_{k,i}
                                                    \|^2
    \end{align}
where the second term in \eqref{Equ:Simulation:L1_MSD_modified} can be evaluated by
expression \eqref{Equ:ConvergenceAnalysis:MSD_Network} with $w^o$ replaced
by $\hat{w}^o$. Moreover, in the derivation of \eqref{Equ:Simulation:L1_MSD_modified},
we used the fact that $\lim_{i\rightarrow \infty}\E(\hat{w}^o-\bm{w}_{k,i})=0$ to eliminate
the cross term, and this result is due to \eqref{Equ:PerformanceAnalysis:Ew_inf_zero} with
$w^o$ there replaced by $\hat{w}^o$.
}%
Fig. \ref{fig:Simulation:LearningCurve} shows the learning curves for different algorithms
for $\gamma=2$ and $\epsilon=10^{-3}$. We see that the diffusion and incremental schemes
have similar performance, and both of them have about $10$ dB gain over the non-cooperation case.
To examine the impact of the parameter $\epsilon$ and the regularization factor $\gamma$,
we show the steady-state MSE for different values of $\gamma$
and $\epsilon$ in Fig. \ref{fig:Simulation:MSDvsGamma}.
When $\epsilon$ is small ($\epsilon=10^{-2}$), adding a reasonable regularization
($\gamma=1\sim4$) decreases the steady-state MSE.  However, when $\epsilon$ is large ($\epsilon=1$),
expression \eqref{Equ:Simulation:R_w} is no longer a good approximation for $\|w\|_1$, and
regularization does not improve the MSE.

\subsection{Distributed Collaborative Localization}
\label{Sec:Simulation:Localization}

The previous example deals with a convex cost \eqref{Equ:Simulation:J_glob_L1RLS}.
Now, we consider a localization problem
that has a non-convex cost function and apply the same diffusion strategies to its solution. Assume each node
is interested in locating a common target
located at $w^o=[0 \; 0]^T$.
Each node $k$ knows its position $x_k$ and has a noisy measurement of the squared distance to
the target:
		\begin{align}
            \label{Equ:Simulation:Localization_datamodel}
			\bm{d}_{k}(i)	=	\|w^o - x_k\|^2 + \bm{z}_k(i), \quad k=1,2,\ldots,N	\nonumber
		\end{align}
where $\bm{z}_k(i) \sim \mathcal{N}(0,\sigma_{z,k}^2)$ is the measurement noise of
node $k$ at time $i$.
The
component cost function $J_k(w)$ at node $k$ is chosen as
	\begin{align}
		J_k(w)		=	\frac{1}{4}\mathbb{E}\left|\bm{d}_k(i) - \|w - x_k\|^2 \right|^2
	\end{align}
where we multiply by $1/4$ here to eliminate a factor of $4$ that will otherwise appear in the
gradient.
If each node $k$ minimizes $J_k(w)$ individually,
it is not possible to solve for $w^o$.
Therefore, we use
information from other nodes, and instead
seek to minimize the following global
cost:
	\begin{align}
		\label{Equ:Simulation:Example1:J_glob}
		J^{\mathrm{glob}}(w)	=	\frac{1}{4}
                                    \sum_{k=1}^N \mathbb{E}\left|\bm{d}_k(i) - \|w - x_k\|^2 \right|^2
	\end{align}
This problem arises, for example, in cellular communication systems, where
multiple base-stations are interested in locating users using the measured
distances between themselves and the user.
Diffusion algorithms \eqref{Equ:DiffusionAdaptation:ATC0}
and \eqref{Equ:DiffusionAdaptation:CTA0}
can be applied to solve the problem in a distributed manner.
Each node $k$ would update its estimate of ${w}^o$ by using
the gradient vectors of $\{J_l(w)\}_{l \in \mathcal{N}_k}$, which are given by:
	\begin{align}
		\nabla_w J_l(w)		=	-\mathbb{E}\bm{d}_l(i) \; (w-x_l)
							+\|w-x_l\|^2 (w-x_l)
	\end{align}
However, the nodes are assumed to have access to measurements $\{{d}_{l}(i),x_l\}$
and not to $\mathbb{E}\bm{d}_l(i)$. In this case, nodes can use the available
measurements to approximate the gradient vectors in
\eqref{Equ:DiffusionAdaptation:ATC}
and \eqref{Equ:DiffusionAdaptation:CTA} as:
	\begin{align}
			\widehat{\nabla}_w J_l(w)
				=	-{d}_l(i) (w-x_l)
							+ \|w-x_l\|^2 (w-x_l)
	\end{align}
If we do not exchange the local gradients with neighbors, i.e., if we set $S=I$,
then the base-stations only share the local estimates of the target position $w^o$
with their neighbors (no exchange of $\{x_l\}_{l \in \mathcal{N}_k}$).

\begin{figure*}[t!]
   \centerline
   {
   	\subfigure[Learning curves for stationary target ($\mu=0.0025$).]
	{
	   \includegraphics[width=3in]{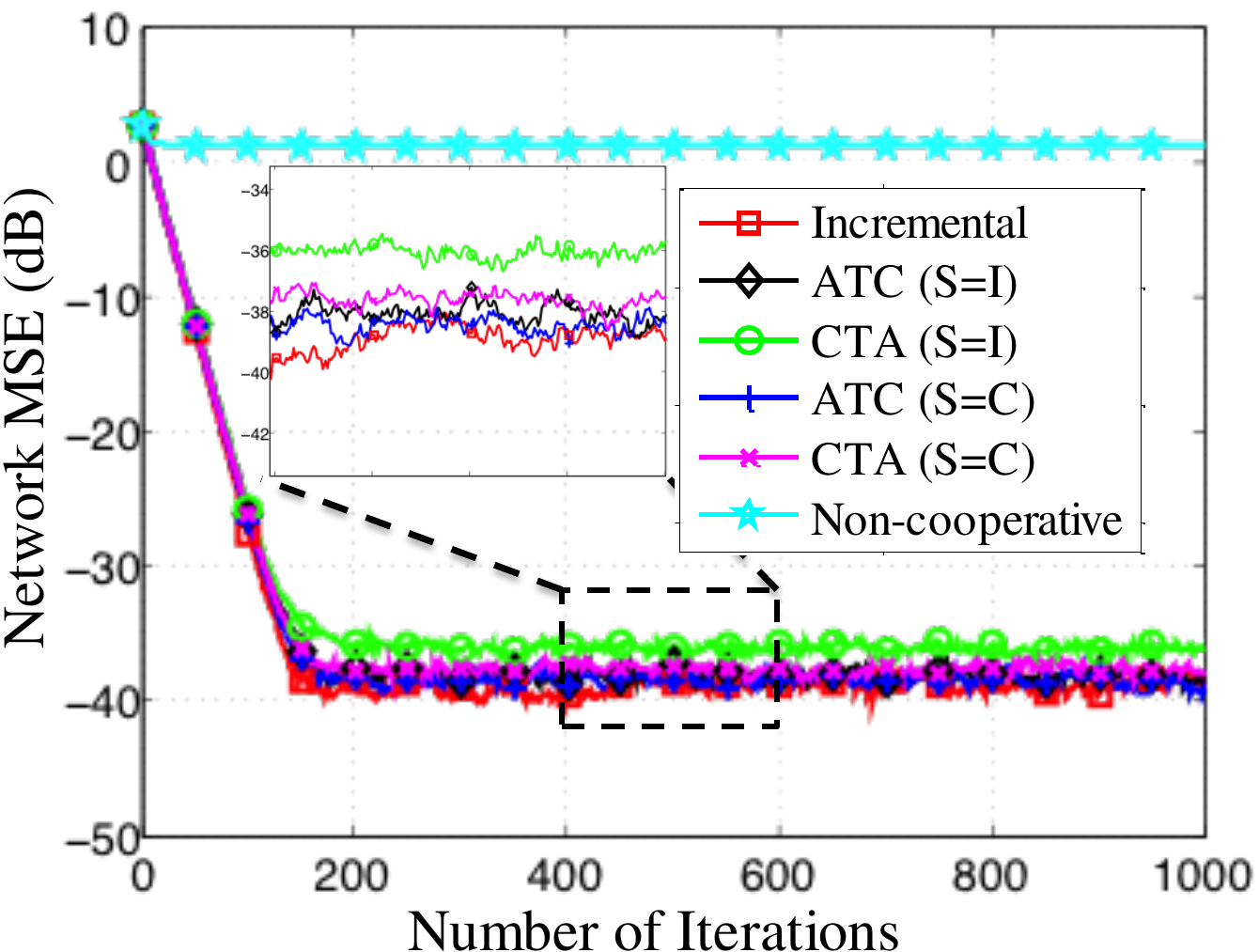}
	   \label{fig:Fig_Localization_LearningCurve_MSD_Truewo}
	}
    \hfil
	\subfigure[Steady-state performance for stationary target.]
	{	   \includegraphics[width=2.9in]{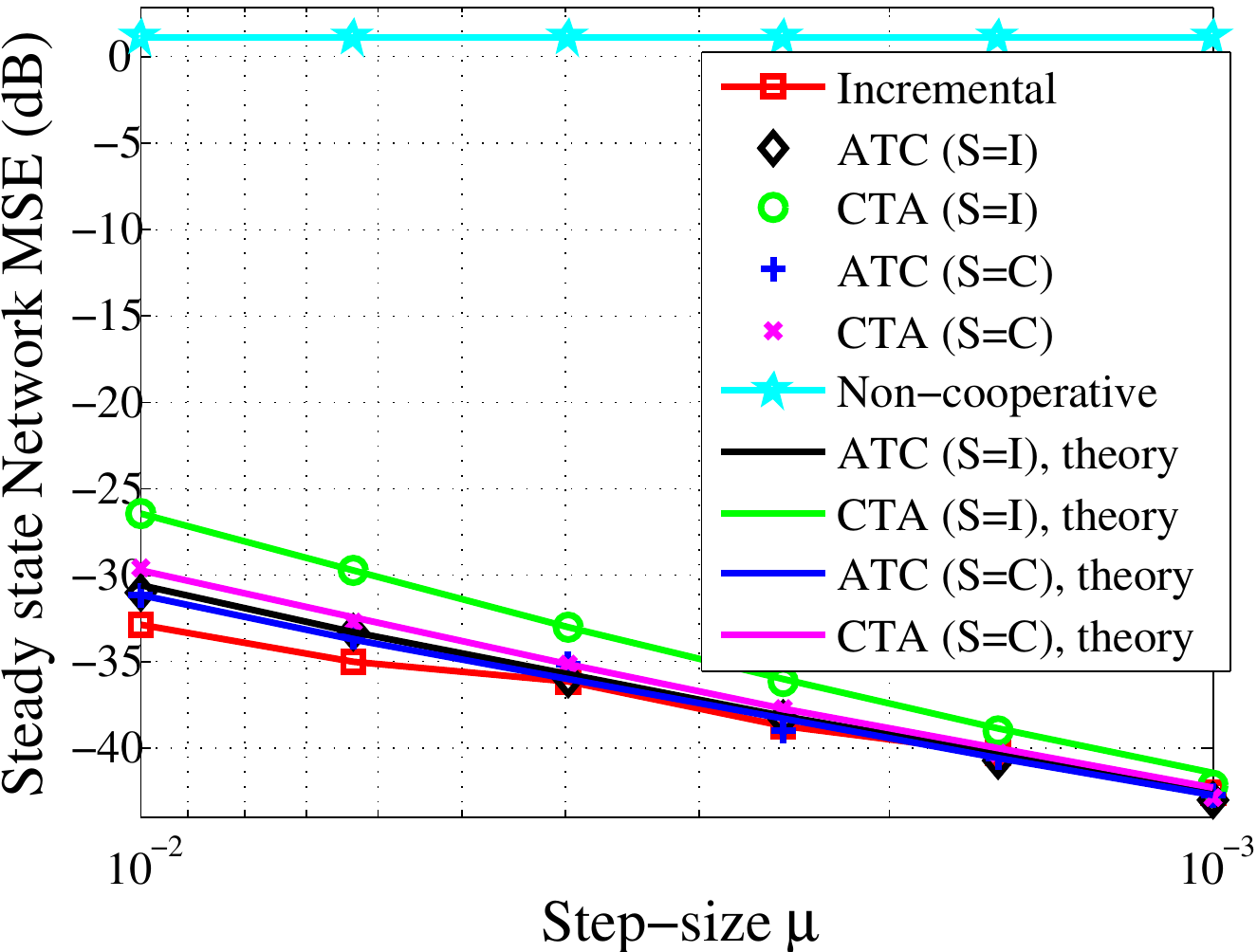}
	   \label{fig:Fig_Localization_MSDvsMU_Truewo}
	  }	
	}
   \centerline
   {
   	\subfigure[Tracking a moving-target by node $1$ ($\mu=0.01$).]
	{
	   \includegraphics[width=3in]{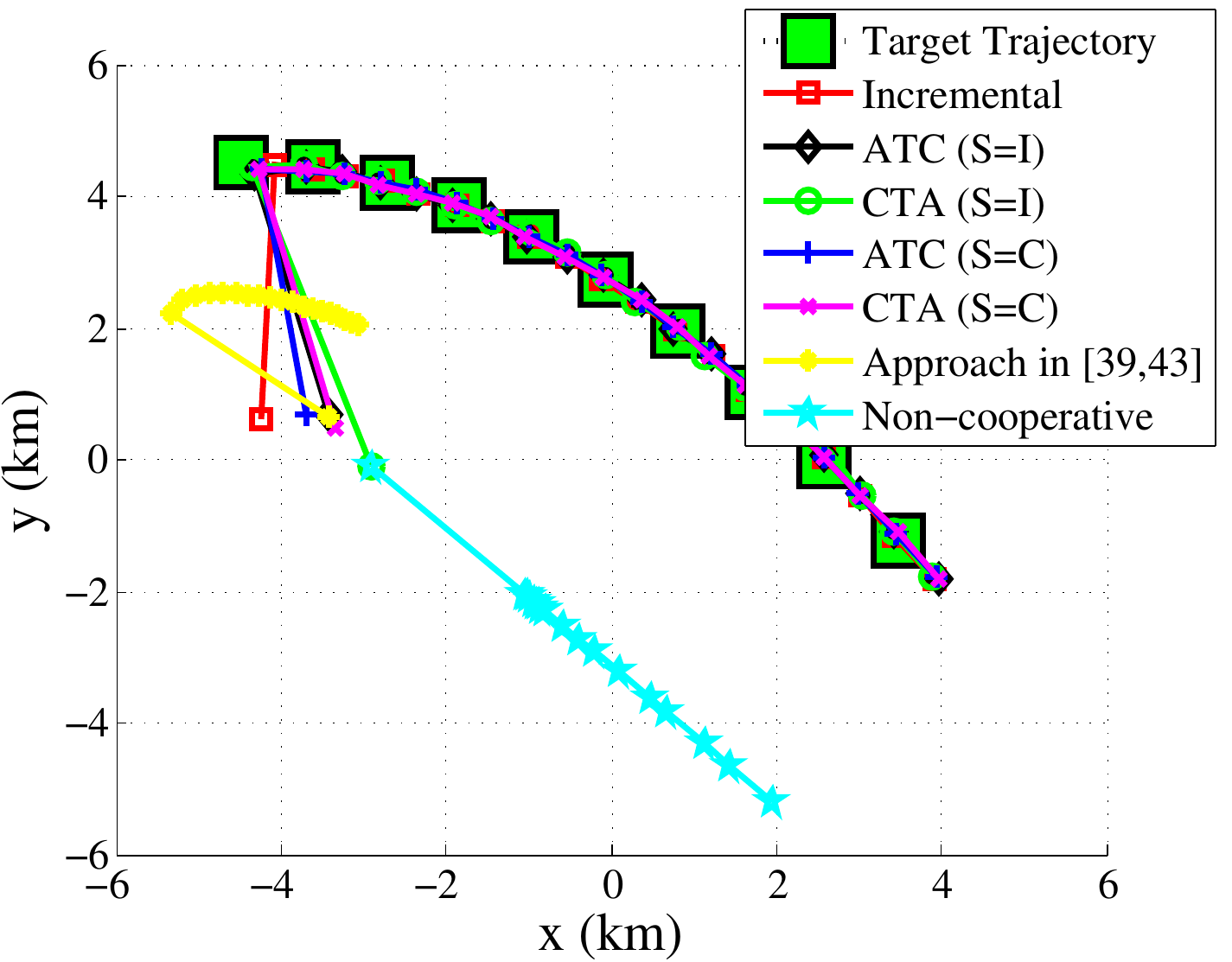}
	   \label{fig:Fig_Localization_TrackingTrajectory}
   	}
    \hfil
	\subfigure[Learning curves for moving target ($\mu=0.01$).]
	{
	   \includegraphics[width=3in]{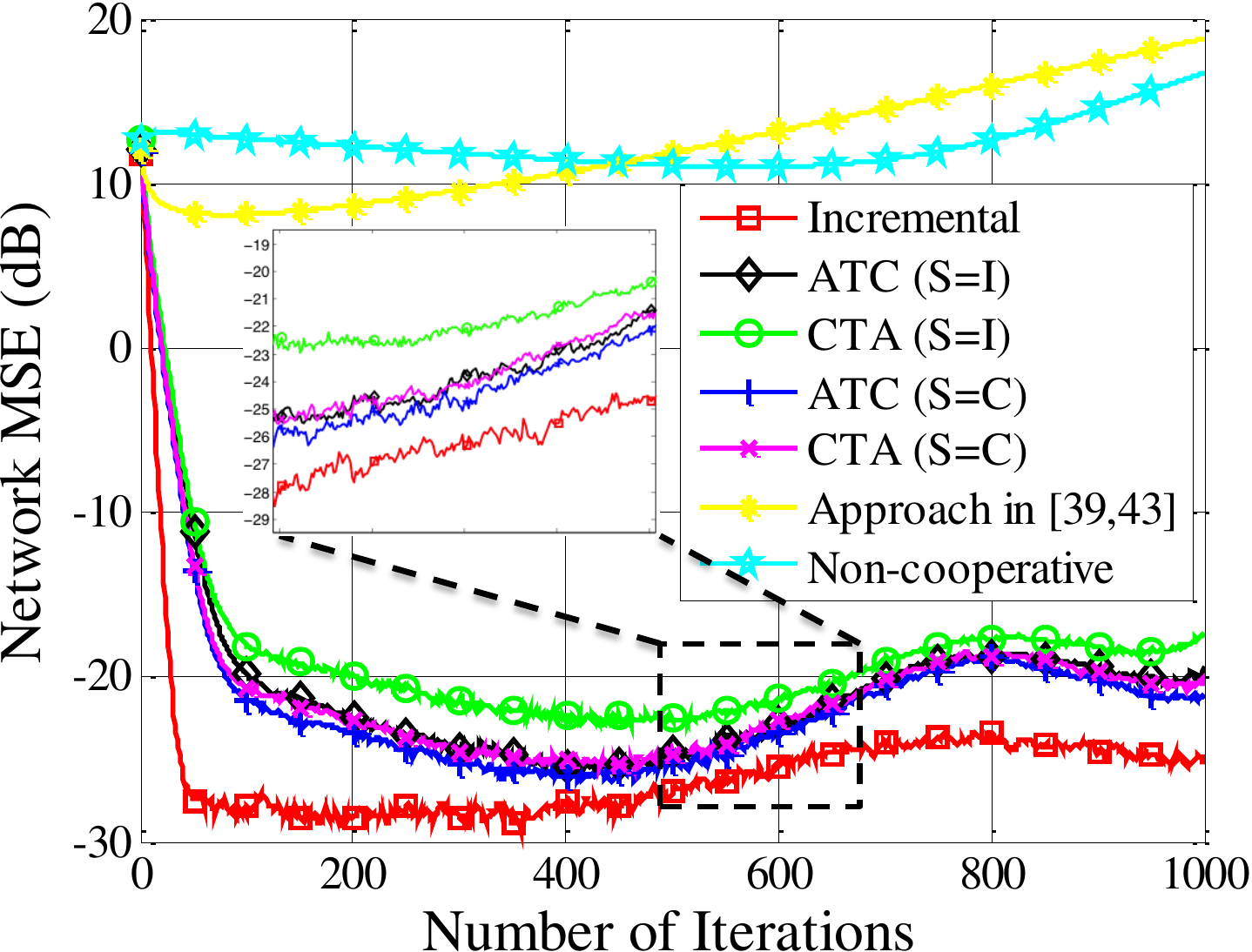}
	   \label{fig:Fig_Localization_TrackingMSD}
   	}
   }
   \caption{Performance of distributed localization for stationary and moving targets.
   	Diffusion strategies employ constant step-sizes, which enable continuous
	   adaptation and learning even when the target moves
	   (which corresponds to a changing cost function).}
\end{figure*}

We first simulate the stationary case, where the target stays at $w^o$.
In Fig. \ref{fig:Fig_Localization_LearningCurve_MSD_Truewo}, we show the
MSE curves for non-cooperative, ATC, CTA, and incremental
algorithms. The noise variance is set to $\sigma_{z,k}^2=1$. We set the step-sizes
to $\mu=0.0025/N$ for the incremental algorithm, and $\mu= 0.0025$ for other algorithms.
For ATC and CTA strategies, we use simple
averaging for the combination step $\{a_{l,k}\}$, and for ATC and CTA with
gradient exchange, we use Metropolis weights for $\{c_{l,k}\}$ to combine the gradients.
The performance of CTA and ATC algorithms are close
to each other, and both of them are
close to the incremental scheme. In Fig. \ref{fig:Fig_Localization_MSDvsMU_Truewo},
we show the steady state MSE with respect to different values of $\mu$.
{\color{black}
We also use expression \eqref{Equ:ConvergenceAnalysis:MSD_Network}
to evaluate the theoretical performance of the diffusion strategies.
}
As the step-size becomes small, the performances of diffusion and incremental algorithms are close,
and the MSE decreases as $\mu$ decreases.
Furthermore, we see that exchanging only local estimates ($S=I$) is enough for
localization, compared to the case of exchanging both local estimates and gradients ($S=C$).

Next, we apply the algorithms to a non-stationary scenario, where the target moves along a trajectory,
as shown in Fig. \ref{fig:Fig_Localization_TrackingTrajectory}. The step-size is set to $\mu=0.01$
for diffusion algorithms, and to $\mu=0.01/N$ for the incremental approach. To see the advantage
of using a constant step-size for continuous tracking, we also simulate the vanishing
step-size version of the algorithm from \cite{ram2010distributed,srivastava2011distributed} ($\mu_{k,i} = 0.01/i$).
The diffusion algorithms track well the target
but not the non-cooperative algorithm and the algorithm from
\cite{ram2010distributed,srivastava2011distributed},
because a decaying step-size is not helpful for tracking.
The tracking performance
is shown in Fig. \ref{fig:Fig_Localization_TrackingMSD}.

\section{Conclusion}
\label{Sec:Conclusion}
This paper proposed diffusion adaptation strategies to optimize global cost functions
over a network of nodes, where the cost consists of several components.
Diffusion adaptation allows the nodes to solve the distributed
optimization problem via local interaction and online learning.
We used gradient approximations and constant step-sizes to endow the networks with continuous
learning and tracking abilities. We analyzed the mean-square-error performance of the algorithms
in some detail,
including their transient and steady-state behavior. Finally,
we applied the scheme to two examples: distributed sparse parameter estimation
and distributed localization. Compared to incremental methods, diffusion strategies do
not require a cyclic path over the nodes, which makes them more robust to node and link failure.

\appendices

\section{Proof of Mean-Square Stability}
\label{Appendix:Proof_MeanSquaredStability}
Taking the $\infty-$norm of both sides of
\eqref{Equ:ConvergenceAnalysis:VarianceRecursion_Inequality_Final}, we obtain
{
    \begin{align}
        \|\mathcal{W}_i\|_{\infty}  \;\le&\;
                                        \|P_2^T \Gamma P_1^T\|_{\infty}
                                        \!\cdot\!
                                        \|\mathcal{W}_{i-1}\|_{\infty}
                                        \!+\!
                                        \sigma_v^2
                                        \|S\|_1^2
                                        \!\cdot\!
                                        \|P_2^T\|_{\infty}
                                        \!\cdot\!
                                        \|\Omega\|_{\infty}^2
                                        \nonumber\\
                                    \;\le&\;
                                        \|P_2^T\|_{\infty}
                                        \cdot
                                        \|\Gamma \|_{\infty}
                                        \cdot
                                        \|P_1^T\|_{\infty}
                                        \cdot
                                        \|\mathcal{W}_{i-1}\|_{\infty}
                                        \nonumber\\
                                        &\;+
                                        \sigma_v^2
                                        \|S\|_1^2
                                        \cdot
                                        \|P_2^T\|_{\infty}
                                        \cdot
                                        \|\Omega\|_{\infty}^2
                                        \nonumber\\
        \label{Equ:Appendix:W_i_norm_inf_ineq}
                                    \;=&\;
                                        \|\Gamma \|_{\infty}
                                        \!\cdot\!
                                        \|\mathcal{W}_{i-1}\|_{\infty}
                                        \!+\!
                                        \big(\max_{1 \le k \le N} \mu_k^2\big)
                                        \!\cdot\!
                                        \sigma_v^2
                                        \|S\|_1^2
    \end{align}
}%
where we used the fact that $\|P_1^T\|_{\infty}=\|P_2^T\|_{\infty}=1$ because
each row of $P_1^T$ and $P_2^T$ sums up to one. Moreover,
from \eqref{Equ:ConvergenceAnalysis:Gamma}, we have
			\begin{align}
				\label{Equ:ConvergenceAnalysis:gamma}
				\|\Gamma\|_{\infty}	=	\max_{1 \le k \le N} (\gamma_k^2+\mu_k^2\alpha\|S\|_1^2)
			\end{align}
Iterating \eqref{Equ:Appendix:W_i_norm_inf_ineq}, we obtain
    \begin{align}
        \|\mathcal{W}_i\|_{\infty}  \;\le\;&
                                        \|\Gamma\|_{\infty}^i \cdot \|\mathcal{W}_0\|_{\infty}
                                        \nonumber\\
        \label{Equ:ConvergenceAnalysis:TheoremProof:W_i_Expr}
                                        \;&+
                                        \big(\max_{1 \le k \le N} \mu_k^2\big)
                                        \cdot
                                        \sigma_v^2
                                        \|S\|_1^2
                                        \sum_{j=0}^{i-1} \|\Gamma\|_{\infty}^j
    \end{align}
We are going to show further ahead that condition
\eqref{Equ:ConvergenceAnalysis:MS_Theorem:Stepsizes}
guarantees $\|\Gamma\|_{\infty}<1$. Now, given that $\|\Gamma\|_{\infty}<1$,
the first term on the right hand side of \eqref{Equ:ConvergenceAnalysis:TheoremProof:W_i_Expr}
converges to zero as $i \rightarrow \infty$, and
the second term	on the right-hand side of \eqref{Equ:ConvergenceAnalysis:TheoremProof:W_i_Expr}
converges to:
			\begin{align}
				\label{Equ:ConvergenceAnalysis:Convergence_SecondTerm_Entry_k}
				\lim_{i \rightarrow \infty}
                                        \sigma_v^2
                                        \|S\|_1^2
                                        \sum_{j=0}^{i-1} \|\Gamma\|_{\infty}^j
                                        =
                                            \frac{\sigma_v^2\|S\|_1^2}{1-\|\Gamma\|_{\infty}}
			\end{align}
Therefore, we establish \eqref{Equ:ConvergenceAnalysis:MS_Theorem:W_infty_Bound} as follows:
			\begin{align}
				\limsup_{i\rightarrow \infty}
                \|\mathcal{W}_{i}\|_{\infty}
                                        &\le
                                            \big(\max_{1 \le k \le N} \mu_k^2\big)
                                            \cdot
                                            \frac{\sigma_v^2\|S\|_1^2}{1-\|\Gamma\|_{\infty}}
                                            \nonumber\\
                \label{Equ:Appendix:lisup_W_i_inf_norm_bound}
										&=	\frac{
													\displaystyle
													\Big(\max_{1\le k \le N} \mu_k^2\Big)
                                                    \cdot
                                                    \|S\|_1^2\sigma_v^2
												}
												{
													\displaystyle
													1
													-
													\max_{1 \le k \le N}
													(\gamma_k^2 + \mu_k^2 \alpha\|S\|_1^2)
												}											
			\end{align}
			
		The only fact that remains to prove is to show that
		\eqref{Equ:ConvergenceAnalysis:MS_Theorem:Stepsizes}
		ensures $\|\Gamma\|_{\infty}<1$. From \eqref{Equ:ConvergenceAnalysis:gamma},
		we see that the condition
		$\|\Gamma\|_{\infty}<1$ is equivalent to requiring:
			\begin{align}
				\label{Equ:ConvergenceAnalysis:gamma_bound1}
				&
				\gamma_k^2 + \mu_k^2 \alpha\|S\|_1^2 < 1, \qquad k=1,\ldots,N.
			\end{align}
		Then, using \eqref{Equ:ConvergenceAnalysis:gamma_k}, this is equivalent to:
			\begin{align}
				\label{Equ:Appendix:Condition_Stepsize_intermediate_original1}
				&\big(1-\mu_k \sum_{l=1}^N s_{l,k} \lambda_{l,\max}\big)^2 + \mu_k^2 \alpha\|S\|_1^2 < 1
				\qquad
				\\
				\label{Equ:Appendix:Condition_Stepsize_intermediate_original2}
				&
				\big(1-\mu_k \sum_{l=1}^N s_{l,k} \lambda_{l,\min}\big)^2 + \mu_k^2 \alpha \|S\|_1^2< 1
			\end{align}
		for $k=1,\ldots,N$.
		Recalling the definitions for $\sigma_{k,\max}$ and
		$\sigma_{k,\min}$ in \eqref{Equ:ConvergenceAnalysis:MS_Theorem:sigma_min_max}
		and solving these two quadratic inequalities with respect to $\mu_k$, we
		arrive at:
			\begin{align}				
				&0 < \mu_k < \frac{2\sigma_{k,\max}}{\sigma_{k,\max}^2 + \alpha\|S\|_1^2},
                \qquad
                0 < \mu_k < \frac{2\sigma_{k,\min}}{\sigma_{k,\min}^2 + \alpha\|S\|_1^2}
				\nonumber
			\end{align}
		and we are led to \eqref{Equ:ConvergenceAnalysis:MS_Theorem:Stepsizes}.

\section{Block Maximum Norm of a Matrix}
\label{Appendix:Proof_Lemma_BMNorm_BlockDiagonalMatrix}
Consider a block matrix $X$ with blocks of size $M\times M$ each.
Its block maximum norm is defined as\cite{takahashi2010diffusion}:
	\begin{align}
		\label{Equ:Appendix:BMNorm:Def_BMNorm_Matrix}
		\|X\|_{b,\infty}	\triangleq	\max_{x \neq 0} \frac{\|X x\|_{b,\infty}}{\|x\|_{b,\infty}}
	\end{align}
where the block maximum norm of a vector $x \triangleq \mathrm{col}\{x_1,\ldots,x_N\}$, formed
by stacking $N$ vectors of size $M$ each on top of each other,
is defined as\cite{takahashi2010diffusion}:
	\begin{align}
		\label{Equ:Appendix:BMNorm:Def_BMNorm_Vector}
		\|x\|_{b,\infty}	\triangleq	\max_{1\le k \le N} \|x_k\|
	\end{align}
where $\|\cdot\|$ denotes the Euclidean norm of its vector argument.
	\begin{lemma}[Block maximum norm]
		\label{Lemma:BMNorm_BlockDiagonalMatrix}
		If a block diagonal matrix
        $X \triangleq \mathrm{diag}\{X_1,\ldots,X_N\}\in \mathbb{R}^{NM \times NM}$
		consists of $N$
		blocks along the diagonal with dimension $M \times M$ each,
		then the block maximum norm of $X$ is bounded as
			\begin{align}
				\label{Equ:ConvergenceAnalysis:Lemma:BMNorm_BlockDiagonalMatrix_Bound}
				\|X\|_{b,\infty}	\le	\max_{1\le k \le N} \|X_k\|
			\end{align}
		in terms of the $2$-induced norms of $\{X_k\}$ (largest singular values).
		Moreover, if $X$ is symmetric,
		then equality holds in \eqref{Equ:ConvergenceAnalysis:Lemma:BMNorm_BlockDiagonalMatrix_Bound}.
	\end{lemma}
	\begin{IEEEproof}
        Note that $X  x	\!=\!	\mathrm{col}\{X_1 x_1,\!\ldots,\!X_Nx_N\}$.
        Evaluating the block maximum norm of vector $X x$ leads to
        	\begin{align}
        		\|X x\|_{b,\infty}	&=	\max_{1 \le k \le N} \| X_k x_k \|		
                                        \nonumber\\
        						&\le	\max_{1 \le k \le N} \| X_k \| \cdot \|x_k \|		
                                        \nonumber\\
        		\label{Equ:Appendix:BMNorm:BMNorm_Xx}
        						&\le	\max_{1 \le k \le N} \| X_k \| \cdot \max_{1 \le k \le N}\|x_k \|
        	\end{align}
        Substituting \eqref{Equ:Appendix:BMNorm:BMNorm_Xx} and
        \eqref{Equ:Appendix:BMNorm:Def_BMNorm_Vector} into
        \eqref{Equ:Appendix:BMNorm:Def_BMNorm_Matrix}, we establish
        \eqref{Equ:ConvergenceAnalysis:Lemma:BMNorm_BlockDiagonalMatrix_Bound} as
        	\begin{align}
        		\|X\|_{b,\infty}	&\triangleq	\max_{x \neq 0} \frac{\|X x\|_{b,\infty}}{\|x\|_{b,\infty}}
        								\nonumber\\
        					&\le			\max_{x \neq 0}
        								\frac{\max_{1 \le k \le N} \| X_k \| \cdot \max_{1 \le k \le N}\|x_k \|}
        									{\max_{1\le k \le N} \|x_k\|}
        								\nonumber\\
        		\label{Equ:Appendix:BMNorm:BMNorm_X_bound}
        					&=			\max_{1 \le k \le N} \| X_k \|
        	\end{align}

        Next, we prove that, if all the diagonal blocks of $X$ are symmetric, then
        equality should hold in \eqref{Equ:Appendix:BMNorm:BMNorm_X_bound}.
        To do this, we only need to show that there exists an $x_0 \neq 0$, such that
        	\begin{align}
        		\label{Equ:Appendix:BMNorm:BMNorm_X_attainablepoint}
        		\frac{\|X x_0\|_{b,\infty}}{\|x_0\|_{b,\infty}}
        					=		\max_{1 \le k \le N} \| X_k \|
        	\end{align}
        which would mean that
        	\begin{align}
        		\|X\|_{b,\infty}	&\triangleq	\max_{x \neq 0} \frac{\|X x\|_{b,\infty}}{\|x\|_{b,\infty}}
        								\nonumber\\
        					&\ge		\frac{\|X x_0\|_{b,\infty}}{\|x_0\|_{b,\infty}}
        							\nonumber\\
        		\label{Equ:Appendix:BMNorm:BMNorm_X_lowerbound}
        					&=		\max_{1 \le k \le N} \| X_k \|	
        	\end{align}
        Then, combining inequalities \eqref{Equ:Appendix:BMNorm:BMNorm_X_bound} and
        \eqref{Equ:Appendix:BMNorm:BMNorm_X_lowerbound}, we would obtain desired equality
        that
        	\begin{align}
        		\label{Equ:Appendix:X_Xb_bmaxNorm}
        		\|X\|_{b,\infty}	=	\max_{1 \le k \le N} \|X_k\|
        	\end{align}
        when $X$ is block diagonal and symmetric.
        Thus, without loss of generality, assume the maximum in \eqref{Equ:Appendix:BMNorm:BMNorm_X_attainablepoint}
        is achieved by $X_1$, i.e.,
        	\begin{align}
        		\displaystyle\max_{1 \le k \le N} \|X_k\| = \|X_1\|
        		\nonumber
        	\end{align}
        For a symmetric $X_k$, its 2-induced norm $\|X_k\|$ (defined as
        the largest singular value of $X_k$) coincides with the spectral radius
        of $X_k$. Let $\lambda_{0}$ denote the eigenvalue of $X_1$ of largest magnitude,
        with the corresponding right eigenvector given by $z_0$. Then,
        	\begin{align}
        		\max_{1 \le k \le N} \|X_k\|	=	|\lambda_0|,
        		\qquad
        		X_1 z_0 = \lambda_0 z_0	\nonumber
        	\end{align}
        We select $x_0=\mathrm{col}\{z_0,0,\ldots,0\}$.
        Then, we establish  \eqref{Equ:Appendix:BMNorm:BMNorm_X_attainablepoint} by:
        	\begin{align}
        		\frac{\|X x_0\|_{b,\infty}}{\|x_0\|_{b,\infty}}
        					&=		\frac{\|\mathrm{col}\{X_1 z_0,0,\ldots,0\}\|_{b,\infty}}
        								{\|\mathrm{col}\{z_0,0,\ldots,0\}\|_{b,\infty}}
                                    \nonumber\\
        					&=		\frac{\|X_1 z_0\|}
        								{\|z_0\|}
        					=		\frac{\|\lambda_0 z_0\|}
        								{\|z_0\|}
        					=		|\lambda_0|
        					=		\max_{1 \le k \le N} \|X_k\|
        							\nonumber
    	\end{align}
\end{IEEEproof}

\section{Stability of $\mathcal{B}$ and $\mc{F}$}
\label{Appendix:Stability_F}
Recall the definitions of the matrices $\mc{B}$ and $\mc{F}$ from \eqref{Equ:PerformanceAnalysis:B_cal_def}
and \eqref{Equ:ConvergenceAnalysis:F}:
    \begin{align}
        \label{Equ:ConvergenceAnalysis:SpectralRadius_B}
        \mc{B}  &=      \mc{P}_2^T [I_{MN} - \mc{M}\mc{D}_{\infty}]\mc{P}_1^T       \\
        \mc{F}  &=      \big(\mc{P}_1 [I_{MN} - \mc{M}\mc{D}_{\infty}]\mc{P}_2\big)
                        \otimes
                        \big(\mc{P}_1 [I_{MN} - \mc{M}\mc{D}_{\infty}]\mc{P}_2\big) \nonumber\\
        \label{Equ:ConvergenceAnalysis:SpectralRadius_F}
                &=       \mc{B}^T \otimes \mc{B}^T
    \end{align}
From \eqref{Equ:ConvergenceAnalysis:SpectralRadius_B}--\eqref{Equ:ConvergenceAnalysis:SpectralRadius_F},
we obtain (see Theorem 13.12 from \cite[p.141]{laub2005matrix}):
    \begin{align}
        \rho(\mc{F})    =   \rho(\mc{B}^T \otimes \mc{B}^T)
                        =   [\rho(\mc{B}^T)]^2
                        =   [\rho(\mc{B})]^2
    \end{align}
where $\rho(\cdot)$ denotes the spectral radius of its matrix argument.
Therefore, the stability of the matrix $\mc{F}$ is equivalent to the
stability of the matrix $\mc{B}$, and we only need to examine the stability of $\mc{B}$.
Now note that the block maximum norm (see the definition in Appendix
\ref{Appendix:Proof_Lemma_BMNorm_BlockDiagonalMatrix}) of the matrix $\mc{B}$ satisfies
    \begin{align}
        \|\mc{B}\|_{b,\infty}   \le     \|I_{MN} - \mc{M}\mc{D}_{\infty}\|_{b,\infty}
    \end{align}
since the block maximum norms of $\mc{P}_1$ and $\mc{P}_2$ are one
(see \cite[p.4801]{takahashi2010diffusion}):
    \begin{align}
        &\left\|\mathcal{P}_1^T\right\|_{b,\infty} = 1,     \qquad
        \left\|\mathcal{P}_2^T\right\|_{b,\infty} = 1
    \end{align}
Moreover, by noting that the spectral radius of a matrix is upper bounded by any
matrix norm (Theorem 5.6.9, \cite[p.297]{horn1990matrix})
and that $I_{MN} - \mc{M}\mc{D}_{\infty}$ is symmetric and block diagonal, we have
    \begin{align}
        \label{Equ:Appendix:rho_B_bound}
        \rho(\mc{B})    \le     \|I_{MN} - \mc{M}\mc{D}_{\infty}\|_{b,\infty}
                        =       \rho(I_{MN} - \mc{M}\mc{D}_{\infty})
    \end{align}
Therefore, the stability of $\mc{B}$ is guaranteed by the stability
of $I_{MN} - \mc{M}\mc{D}_{\infty}$.
Next, we call upon the following lemma to
bound $\left\|I_{MN}\!-\!\mathcal{M}{\mathcal{D}}_{\infty}\right\|_{b,\infty}$.
	\begin{lemma}[Norm of $I_{MN}\!-\!\mathcal{M}\mathcal{D}_{\infty}$]
		\label{Lemma:I_MDinf_Bound}
		It holds that the
        matrix $\mathcal{D}_{\infty}$ defined in \eqref{Equ:ConvergenceAnalysis:D_inf}
		satisfies
			\begin{align}
				\label{Equ:ConvergenceAnalysis:I_MDinf_Bound}
						\left\|
							I_{MN}\!-\!\mathcal{M}{\mathcal{D}}_{\infty}
						\right\|_{b,\infty}		
					\le	\max_{1 \le k \le N} \gamma_k
			\end{align}
		where $\gamma_k$ is defined in \eqref{Equ:ConvergenceAnalysis:gamma_k}.
	\end{lemma}
	\begin{IEEEproof}
		Since $\mathcal{D}_{\infty}$ is block diagonal and symmetric, $I_{MN}-\mathcal{M}\mathcal{D}_{\infty}$
		is also block diagonal with blocks $\{I_{M}\!-\!\mu_k{\mathcal{D}}_{k,\infty}\}$,
		where ${\mathcal{D}}_{k,\infty}$ denotes the $k$th diagonal block of $\mathcal{D}_{\infty}$.
		Then, from \eqref{Equ:ConvergenceAnalysis:Lemma:BMNorm_BlockDiagonalMatrix_Bound}
		in Lemma \ref{Lemma:BMNorm_BlockDiagonalMatrix} in Appendix
        \ref{Appendix:Proof_Lemma_BMNorm_BlockDiagonalMatrix},
		it holds that
			\begin{align}
				\label{Equ:ConvergenceAnalysis:I_MDinf_NoriyukiFact}
				\left\|
					I_{MN}\!-\!\mathcal{M}{\mathcal{D}}_{\infty}
				\right\|_{b,\infty}	
						&=	\max_{1 \le k \le N}
							\left\|
								I_{M}\!-\!\mu_k{\mathcal{D}}_{k,\infty}
							\right\|
			\end{align}		
		By the definition of $\mathcal{D}_{\infty}$ in
		\eqref{Equ:ConvergenceAnalysis:D_inf}, and using condition \eqref{Assumption:StrongConvexity}
		from Assumption \ref{Assumption:Hessian}, we have
			\begin{align}	
				\Big(\sum_{l=1}^N s_{l,k} \lambda_{l,\min}\Big) \cdot I_M
				\le
				{\mathcal{D}}_{k,\infty}	
				\le	
				\Big(\sum_{l=1}^N s_{l,k} \lambda_{l,\max}\Big) \cdot I_M
				\nonumber
			\end{align}
		which implies that
            \begin{align}
                \label{Equ:Appendix:Bound_I_muD_k_Lower}
                &I_M - \mu_k \mc{D}_{k,\infty}
                        \ge     \Big(1-\mu_k\sum_{l=1}^N s_{l,k} \lambda_{l,\max}\Big) \cdot I_M    \\
                \label{Equ:Appendix:Bound_I_muD_k_Upper}
                &I_M - \mu_k \mc{D}_{k,\infty}
                        \le     \Big(1- \mu_k\sum_{l=1}^N s_{l,k} \lambda_{l,\min}\Big) \cdot I_M
            \end{align}
        Thus, $\| I_M \!-\! \mu_k \mathcal{D}_{k,\infty} \|	\!\le\!	\gamma_k$. Substituting
		into \eqref{Equ:ConvergenceAnalysis:I_MDinf_NoriyukiFact}, we get
		\eqref{Equ:ConvergenceAnalysis:I_MDinf_Bound}.
	\end{IEEEproof}
	\vspace{0.5em}	
\noindent
Substituting \eqref{Equ:ConvergenceAnalysis:I_MDinf_Bound} into
\eqref{Equ:Appendix:rho_B_bound}, we get:
	\begin{align}
		\rho(\mc{B})	\le	\max_{1 \le k \le N} \gamma_k
	\end{align}
As long as $\displaystyle\max_{1 \le k \le N} \gamma_k < 1$, then all the eigenvalues
of $\mc{B}$ will lie within the unit circle.
By the definition of $\gamma_k$ in \eqref{Equ:ConvergenceAnalysis:gamma_k},
this is equivalent to requiring
    \begin{align}
		|1-\mu_k \sigma_{k,\max}| < 1, \qquad
		|1-\mu_k \sigma_{k,\min}| < 1
        \nonumber
	\end{align}
for $k=1,\ldots,N$,
where $\sigma_{k,\max}$ and $\sigma_{k,\min}$ are defined in
\eqref{Equ:ConvergenceAnalysis:MS_Theorem:sigma_min_max}.
These conditions are satisfied if we choose $\mu_k$ such that
	\begin{align}
		0	<	\mu_k	<	{2}/{\sigma_{k,\max}}, \qquad k=1,\ldots,N
	\end{align}
which is obviously guaranteed for sufficiently small step-sizes (and also by condition
\eqref{Equ:ConvergenceAnalysis:MS_Theorem:Stepsizes}).



\bibliographystyle{IEEEtran}
\bibliography{DistOpt}

\begin{thebibliography}{10}
\providecommand{\url}[1]{#1}
\csname url@rmstyle\endcsname
\providecommand{\newblock}{\relax}
\providecommand{\bibinfo}[2]{#2}
\providecommand\BIBentrySTDinterwordspacing{\spaceskip=0pt\relax}
\providecommand\BIBentryALTinterwordstretchfactor{4}
\providecommand\BIBentryALTinterwordspacing{\spaceskip=\fontdimen2\font plus
\BIBentryALTinterwordstretchfactor\fontdimen3\font minus
  \fontdimen4\font\relax}
\providecommand\BIBforeignlanguage[2]{{%
\expandafter\ifx\csname l@#1\endcsname\relax
\typeout{** WARNING: IEEEtran.bst: No hyphenation pattern has been}%
\typeout{** loaded for the language `#1'. Using the pattern for}%
\typeout{** the default language instead.}%
\else
\language=\csname l@#1\endcsname
\fi
#2}}

\bibitem{chen2011diffusionOpt}
J.~Chen, S.-Y. Tu, and A.~H. Sayed, ``Distributed optimization via diffusion
  adaptation,'' in \emph{Proc. IEEE International Workshop on Comput. Advances
  Multi-Sensor Adaptive Process. (CAMSAP)}, Puerto Rico, Dec. 2011, pp.
  281--284.

\bibitem{chen2011MSE}
J.~Chen and A.~H. Sayed, ``Performance of diffusion adaptation for
  collaborative optimization,'' in \emph{Proc. IEEE International Conf.
  Acoustics, Speech and Signal Process. (ICASSP)}, Kyoto, Japan, March 2012,
  pp. 1--4.

\bibitem{tu2011fish}
S.-Y. Tu and A.~H. Sayed, ``Mobile adaptive networks,'' \emph{IEEE J. Sel.
  Topics. Signal Process.}, vol.~5, no.~4, pp. 649--664, Aug. 2011.

\bibitem{dekel2011optimal}
O.~Dekel, R.~Gilad-Bachrach, O.~Shamir, and L.~Xiao, ``Optimal distributed
  online prediction,'' in \emph{Proc. International Conf. Machin. Learning
  (ICML)}, Bellevue, USA, June 2011, pp. 713--720.

\bibitem{zaid2011collaborativeGMM}
Z.~J. Towfic, J.~Chen, and A.~H. Sayed, ``Collaborative learning of mixture
  models using diffusion adaptation,'' in \emph{Proc. IEEE Workshop on Mach.
  Learning Signal Process. (MLSP)}, Beijing, China, Sep. 2011, pp. 1--6.

\bibitem{bertsekas1997new}
D.~P. Bertsekas, ``A new class of incremental gradient methods for least
  squares problems,'' \emph{SIAM J. Optim.}, vol.~7, no.~4, pp. 913--926, 1997.

\bibitem{nedic2001incremental}
A.~Nedic and D.~P. Bertsekas, ``Incremental subgradient methods for
  nondifferentiable optimization,'' \emph{SIAM J. Optim.}, vol.~12, no.~1, pp.
  109--138, 2001.

\bibitem{rabbat2005quantized}
M.~G. Rabbat and R.~D. Nowak, ``Quantized incremental algorithms for
  distributed optimization,'' \emph{IEEE J. Sel. Areas Commun.}, vol.~23,
  no.~4, pp. 798--808, Apr. 2005.

\bibitem{lopes2007incremental}
C.~G. Lopes and A.~H. Sayed, ``{Incremental adaptive strategies over
  distributed networks},'' \emph{IEEE Trans. Signal Process.}, vol.~55, no.~8,
  pp. 4064--4077, Aug. 2007.

\bibitem{bertsekasparallel}
D.~P. Bertsekas and J.~N. Tsitsiklis, \emph{{Parallel and Distributed
  Computation: Numerical Methods, 1st edition}}.\hskip 1em plus 0.5em minus
  0.4em\relax Athena Scientific, Singapore, 1997.

\bibitem{tsitsiklis1984convergence}
J.~N. Tsitsiklis and M.~Athans, ``Convergence and asymptotic agreement in
  distributed decision problems,'' \emph{{IEEE Trans. Autom. Control}},
  vol.~29, no.~1, pp. 42--50, Jan. 1984.

\bibitem{tsitsiklis1986distributed}
J.~N. Tsitsiklis, D.~P. Bertsekas, and M.~Athans, ``Distributed asynchronous
  deterministic and stochastic gradient optimization algorithms,'' \emph{IEEE
  Trans. Autom. Control}, vol.~31, no.~9, pp. 803--812, Sep. 1986.

\bibitem{barbarossa2007bio}
S.~Barbarossa and G.~Scutari, ``Bio-inspired sensor network design,''
  \emph{IEEE Signal Process. Mag.}, vol.~24, no.~3, pp. 26--35, May 2007.

\bibitem{nedic2009bookchapter}
A.~Nedic and A.~Ozdaglar, ``Cooperative distributed multi-agent optimization,''
  \emph{in Convex Optimization in Signal Processing and Communications, {\rm Y.
  Eldar and D. Palomar, Eds.}}, pp. 340–--386, 2009.

\bibitem{nedic2009distributed}
------, ``Distributed subgradient methods for multi-agent optimization,''
  \emph{IEEE Trans. Autom. Control}, vol.~54, no.~1, pp. 48--61, Jan. 2009.

\bibitem{schizas2008consensus1}
I.~D. Schizas, A.~Ribeiro, and G.~B. Giannakis, ``{Consensus in ad hoc WSNs
  with noisy links---Part I: Distributed estimation of deterministic
  signals},'' \emph{IEEE Trans. Signal Process.}, vol.~56, no.~1, pp. 350--364,
  2008.

\bibitem{kar2008sensor}
S.~Kar and J.~M.~F. Moura, ``Sensor networks with random links: Topology design
  for distributed consensus,'' \emph{IEEE Trans. Signal Process.}, vol.~56,
  no.~7, pp. 3315--3326, July 2008.

\bibitem{kar2011converegence}
------, ``Convergence rate analysis of distributed gossip (linear parameter)
  estimation: Fundamental limits and tradeoffs,'' \emph{IEEE J. Sel. Topics.
  Signal Process.}, vol.~5, no.~4, pp. 674--690, Aug. 2011.

\bibitem{dimakis2010gossip}
A.~G. Dimakis, S.~Kar, J.~M.~F. Moura, M.~G. Rabbat, and A.~Scaglione, ``Gossip
  algorithms for distributed signal processing,'' \emph{Proc. of the IEEE},
  vol.~98, no.~11, pp. 1847--1864, 2010.

\bibitem{olfati2004consensus}
R.~Olfati-Saber and R.~M. Murray, ``Consensus problems in networks of agents
  with switching topology and time-delays,'' \emph{IEEE Trans. Autom. Control},
  vol.~49, no.~9, pp. 1520--1533, Sep. 2004.

\bibitem{aysal2009broadcast}
T.~C. Aysal, M.~E. Yildiz, A.~D. Sarwate, and A.~Scaglione, ``Broadcast gossip
  algorithms for consensus,'' \emph{IEEE Trans. Signal Process.}, vol.~57,
  no.~7, pp. 2748--2761, 2009.

\bibitem{sardellitti2010fast}
S.~Sardellitti, M.~Giona, and S.~Barbarossa, ``Fast distributed average
  consensus algorithms based on advection-diffusion processes,'' \emph{IEEE
  Trans. Signal Process.}, vol.~58, no.~2, pp. 826--842, Feb. 2010.

\bibitem{xiao2005scheme}
L.~Xiao, S.~Boyd, and S.~Lall, ``A scheme for robust distributed sensor fusion
  based on average consensus,'' in \emph{Proc. Int. Symp. Information
  Processing Sensor Networks (IPSN)}, Los Angeles, CA, Apr. 2005, pp. 63--70.

\bibitem{eksin2011asilomar}
C.~Eksin and A.~Ribeiro, ``Network optimization with heuristic rational
  agents,'' in \emph{Proc. Asilomar Conf. on Signals Systems Computers},
  Pacific Grove, CA, Nov. 2011, pp. 1--5.

\bibitem{karp1972reducibility}
R.~M. Karp, ``Reducibility among combinational problems,'' \emph{Complexity of
  Computer Computations {\rm(R. E. Miller and J. W. Thatcher, Eds.)}}, pp.
  85--104, 1972.

\bibitem{lopesdistributed}
C.~G. Lopes and A.~H. Sayed, ``Distributed processing over adaptive networks,''
  in \emph{Proc. Adaptive Sensor Array Processing Workshop}, MIT Lincoln
  Laboratory, MA, June 2006, pp. 1--5.

\bibitem{lopes2007diffusion}
C.~Lopes and A.~Sayed, ``Diffusion least-mean squares over adaptive networks,''
  in \emph{IEEE ICASSP}, vol.~3, Honolulu, HI, Apr. 2007, pp. 917--920.

\bibitem{Sayed07}
A.~H. Sayed and C.~G. Lopes, ``Adaptive processing over distributed networks,''
  \emph{IEICE Trans. Fund. Electron., Commun. Comput. Sci.}, vol. E90-A, no.~8,
  pp. 1504--1510, Aug. 2007.

\bibitem{lopes2008diffusion}
C.~G. Lopes and A.~H. Sayed, ``{Diffusion least-mean squares over adaptive
  networks: Formulation and performance analysis},'' \emph{IEEE Trans. Signal
  Process.}, vol.~56, no.~7, pp. 3122--3136, July 2008.

\bibitem{cattivelli2008diffusion}
F.~S. Cattivelli and A.~H. Sayed, ``{Diffusion LMS algorithms with information
  exchange},'' in \emph{Proc. Asilomar Conf. Signals, Syst. Comput.}, Pacific
  Grove, CA, Nov. 2008, pp. 251--255.

\bibitem{Cattivelli10}
------, ``{Diffusion LMS strategies for distributed estimation},'' \emph{IEEE
  Trans. Signal Process.}, vol.~58, no.~3, pp. 1035--1048, March 2010.

\bibitem{cattivelli2007diffusionRLS}
F.~S. Cattivelli, C.~G. Lopes, and A.~H. Sayed, ``{A diffusion RLS scheme for
  distributed estimation over adaptive networks},'' in \emph{Proc. IEEE
  Workshop on Signal Process. Advances Wireless Comm. (SPAWC)}, Helsinki,
  Finland, June 2007, pp. 1--5.

\bibitem{cattivelli2008TSPdiffusionRLS}
------, ``Diffusion recursive least-squares for distributed estimation over
  adaptive networks,'' \emph{IEEE Trans. Signal Process.}, vol.~56, no.~5, pp.
  1865--1877, May 2008.

\bibitem{cattivelli2010TACdiffusionKalman}
F.~S. Cattivelli and A.~H. Sayed, ``Diffusion strategies for distributed
  {Kalman} filtering and smoothing,'' \emph{IEEE Trans. Autom. Control},
  vol.~55, no.~9, pp. 2069--2084, Sep. 2010.

\bibitem{takahashi2010diffusion}
N.~Takahashi, I.~Yamada, and A.~H. Sayed, ``{Diffusion least-mean squares with
  adaptive combiners: Formulation and performance analysis},'' \emph{IEEE
  Trans. Signal Process.}, vol.~58, no.~9, pp. 4795--4810, Sep. 2010.

\bibitem{cattivelli2011modeling}
F.~S. Cattivelli and A.~H. Sayed, ``Modeling bird flight formations using
  diffusion adaptation,'' \emph{IEEE Trans. Signal Process.}, vol.~59, no.~5,
  pp. 2038--2051, May 2011.

\bibitem{di2011bio}
P.~Di~Lorenzo, S.~Barbarossa, and A.~H. Sayed, ``Bio-inspired swarming for
  dynamic radio access based on diffusion adaptation,'' in \emph{Proc. European
  Signal Process. Conf. (EUSIPCO)}, Aug. 2011, pp. 1--6.

\bibitem{chouvardas2011adaptive}
S.~Chouvardas, K.~Slavakis, and S.~Theodoridis, ``Adaptive robust distributed
  learning in diffusion sensor networks,'' \emph{IEEE Trans. Signal Process.},
  vol.~59, no.~10, pp. 4692--4707, 2011.

\bibitem{ram2010distributed}
S.~S. Ram, A.~Nedic, and V.~V. Veeravalli, ``Distributed stochastic subgradient
  projection algorithms for convex optimization,'' \emph{J. Optim. Theory
  Appl.}, vol. 147, no.~3, pp. 516--545, 2010.

\bibitem{bianchi2011convergence}
P.~Bianchi, G.~Fort, W.~Hachem, and J.~Jakubowicz, ``Convergence of a
  distributed parameter estimator for sensor networks with local averaging of
  the estimates,'' in \emph{Proc. IEEE ICASSP}, Prague, Czech, May 2011, pp.
  3764--3767.

\bibitem{bertsekas2010incremental}
D.~P. Bertsekas, ``Incremental gradient, subgradient, and proximal methods for
  convex optimization: A survey,'' \emph{LIDS Technical Report, MIT}, no. 2848,
  2010.

\bibitem{borkar2000ode}
V.~S. Borkar and S.~P. Meyn, ``The {ODE} method for convergence of stochastic
  approximation and reinforcement learning,'' \emph{SIAM J. Control Optim.},
  vol.~38, no.~2, pp. 447--469, 2000.

\bibitem{srivastava2011distributed}
K.~Srivastava and A.~Nedic, ``Distributed asynchronous constrained stochastic
  optimization,'' \emph{IEEE J. Sel. Topics. Signal Process.}, vol.~5, no.~4,
  pp. 772--790, Aug. 2011.

\bibitem{poliak1987introduction}
B.~Polyak, \emph{{Introduction to Optimization}}.\hskip 1em plus 0.5em minus
  0.4em\relax Optimization Software, NY, 1987.

\bibitem{chen2012ssp}
J.~Chen and A.~H. Sayed, ``{Distributed Pareto-optimal solutions via diffusion
  adaptation},'' in \emph{Proc. IEEE Statistical Signal Process. Workshop
  (SSP)}, Ann Arbor, MI, Aug. 2012.

\bibitem{Sayed08}
A.~H. Sayed, \emph{Adaptive Filters}.\hskip 1em plus 0.5em minus 0.4em\relax
  Wiley, NJ, 2008.

\bibitem{golub1996matrix}
G.~H. Golub and C.~F. Van~Loan, \emph{{Matrix Computations (3rd
  Edition)}}.\hskip 1em plus 0.5em minus 0.4em\relax Johns Hopkins University
  Press, 1996.

\bibitem{stankovic2011decentralized}
S.~S. Stankovic, M.~S. Stankovic, and D.~M. Stipanovic, ``Decentralized
  parameter estimation by consensus based stochastic approximation,''
  \emph{IEEE Trans. Autom. Control}, vol.~56, no.~3, pp. 531--543, Mar. 2011.

\bibitem{Tu2012diffcons}
S.-Y. Tu and A.~H. Sayed, ``Diffusion networks outperform consensus networks,''
  in \emph{Proc. IEEE Statistical Signal Processing Workshop (SSP)}, Ann Arbor,
  MI, Aug. 2012.

\bibitem{horn1990matrix}
R.~Horn and C.~Johnson, \emph{{Matrix Analysis}}.\hskip 1em plus 0.5em minus
  0.4em\relax Cambridge University Press, 1990.

\bibitem{bertsekas2000gradient}
D.~P. Bertsekas and J.~N. Tsitsiklis, ``Gradient convergence in gradient
  methods with errors,'' \emph{SIAM J. Optim.}, vol.~10, no.~3, pp. 627--642,
  2000.

\bibitem{haykin2002adaptive}
S.~Haykin, \emph{{Adaptive Filter Theory}, {2nd Edition}}.\hskip 1em plus 0.5em
  minus 0.4em\relax Prentice Hall, 2002.

\bibitem{arenas2006plant}
J.~Arenas-Garcia, M.~Martinez-Ramon, A.~Navia-Vazquez, and A.~R.
  Figueiras-Vidal, ``Plant identification via adaptive combination of
  transversal filters,'' \emph{Signal Processing}, vol.~86, no.~9, pp.
  2430--2438, 2006.

\bibitem{silva2008improving}
M.~Silva and V.~Nascimento, ``Improving the tracking capability of adaptive
  filters via convex combination,'' \emph{IEEE Trans. Signal Process.},
  vol.~56, no.~7, pp. 3137--3149, 2008.

\bibitem{theodoridis2011adaptive}
S.~Theodoridis, K.~Slavakis, and I.~Yamada, ``Adaptive learning in a world of
  projections,'' \emph{IEEE Signal Process. Mag.}, vol.~28, no.~1, pp. 97--123,
  Jan. 2011.

\bibitem{boyd2004convex}
S.~Boyd and L.~Vandenberghe, \emph{{Convex Optimization}}.\hskip 1em plus 0.5em
  minus 0.4em\relax Cambridge University Press, 2004.

\bibitem{laub2005matrix}
A.~J. Laub, \emph{{Matrix Analysis for Scientists and Engineers}}.\hskip 1em
  plus 0.5em minus 0.4em\relax Society for Industrial and Applied Mathematics
  (SIAM), PA, 2005.

\bibitem{diLorenzo2012icassp}
P.~Di~Lorenzo, S.~Barbarossa, and A.~H. Sayed, ``Sparse diffusion {LMS} for
  distributed adaptive estimation,'' in \emph{Proc. IEEE ICASSP}, Kyoto, Japan,
  March 2012, pp. 1--4.

\bibitem{tibshirani1996regression}
R.~Tibshirani, ``Regression shrinkage and selection via the lasso,'' \emph{J.
  Royal Statist. Soc. B}, pp. 267--288, 1996.

\bibitem{baraniuk2007compressive}
R.~G. Baraniuk, ``Compressive sensing,'' \emph{IEEE Signal Process. Mag.},
  vol.~24, no.~4, pp. 118--121, Mar. 2007.

\bibitem{candes2008enhancing}
E.~Candes, M.~Wakin, and S.~Boyd, ``Enhancing sparsity by reweighted $\ell_1$
  minimization,'' \emph{J. Fourier Anal. Appl.}, vol.~14, no.~5, pp. 877--905,
  2008.

\bibitem{mateos2010distributed}
G.~Mateos, J.~A. Bazerque, and G.~B. Giannakis, ``Distributed sparse linear
  regression,'' \emph{IEEE Trans. Signal Process.}, vol.~58, no.~10, pp.
  5262--5276, 2010.

\bibitem{kopsinis2011online}
Y.~Kopsinis, K.~Slavakis, and S.~Theodoridis, ``Online sparse system
  identification and signal reconstruction using projections onto weighted
  $\ell_1$ balls,'' \emph{IEEE Trans. Signal Process.}, vol.~59, no.~3, pp.
  936--952, Mar. 2011.

\end{thebibliography}

\end{document}